\documentclass[pdflatex,sn-mathphys-num]{sn-jnl}
\usepackage[T1]{fontenc}
\usepackage[utf8]{inputenc}
\usepackage{amsmath,amssymb,amsthm,mathtools}
\usepackage{booktabs}
\usepackage{placeins}
\usepackage[numbers,sort&compress]{natbib}

\theoremstyle{plain}
\newtheorem{theorem}{Theorem}[section]
\newtheorem{proposition}[theorem]{Proposition}
\newtheorem{lemma}[theorem]{Lemma}
\newtheorem{corollary}[theorem]{Corollary}
\theoremstyle{definition}
\newtheorem{definition}[theorem]{Definition}
\newtheorem{assumption}[theorem]{Assumption}
\newtheorem{example}[theorem]{Example}
\theoremstyle{remark}
\newtheorem{remark}[theorem]{Remark}

\newcommand{\HH}{\mathcal{H}}
\newcommand{\Hi}{H}
\newcommand{\Ob}{\mathcal{O}}
\newcommand{\Sa}{\Psi}
\newcommand{\GG}{\mathcal{G}}
\newcommand{\Syn}{\mathcal{S}}
\newcommand{\D}{\mathcal{D}}

\newcommand{\Lacc}{\mathcal{L}}
\newcommand{\N}{\mathbb{N}}
\newcommand{\Z}{\mathbb{Z}}
\newcommand{\R}{\mathbb{R}}
\newcommand{\E}{\mathbb{E}}
\newcommand{\Prob}{\mathbb{P}}
\newcommand{\ip}[2]{\langle #1 , #2 \rangle}
\newcommand{\norm}[1]{\lVert #1 \rVert}
\newcommand{\tr}{\operatorname{tr}}
\newcommand{\ran}{\operatorname{ran}}
\newcommand{\Cov}{\operatorname{Cov}}

\newcommand{\opnorm}[1]{\norm{#1}_{\mathrm{op}}}

\begin{document}

\title{Augmented Dynamical Sampling under unmodeled Forcing
}

\author*[1,2]{\fnm{Juan I.} \sur{Giribet}}\email{jgiribet@udesa.edu.ar}
\author[1]{\fnm{Serena} \sur{Feldberg}}

\affil[1]{\orgname{Universidad de San Andrés}, \orgaddress{\street{Vito Dumas 284}, \city{Victoria}, \postcode{B1644BID}, \state{Buenos Aires}, \country{Argentina}}}
\affil[2]{\orgname{Consejo Nacional de Investigaciones Científicas y Técnicas - CONICET}, \orgaddress{ \state{Buenos Aires}, \country{Argentina}}}



\maketitle

\begin{abstract}

We study recovery of an initial state from dynamical samples of a trajectory driven by unknown additive forcing. Measurements against both the original sensors and their images under the adjoint nominal evolution yield one-step residuals that cancel the initial state exactly. These residuals separate forcing estimation from initial-state reconstruction. For a centered Gaussian forcing model with fixed covariance parameters, we combine Gaussian conditioning with canonical dual-frame reconstruction and derive an exact unconditional second-moment identity, including the correlation between the forcing-estimation error and the measurement noise. We also obtain high-probability bounds and identify the irreducible forcing uncertainty through a conditional covariance in the geometry induced by the prior. A hierarchy of augmented sensor systems separates the role of sampling depth in forcing visibility from that of the observation horizon in nominal observability. An image-transport experiment and finite-dimensional and periodic-convolution experiments illustrate the residual correction and the separate roles of sampling depth and observation horizon.

\end{abstract}

\section{Introduction}\label{sec:intro}

Let $\Hi$ be a separable complex Hilbert space, let $A\in B(\Hi)$, and let $\{g_j\}_{j\in J}\subset \Hi$ be a countable family of \emph{sensors}. Dynamical sampling asks when a signal $f\in\Hi$ can be recovered stably from
\begin{equation}\label{eq:classical-data}
  \ip{A^nf}{g_j},\qquad j\in J,\ 0\le n\le N .
\end{equation}
Since $\ip{A^nf}{g_j}=\ip{f}{(A^*)^ng_j}$, this is a frame question for the iterated system
\begin{equation}\label{eq:frame-system}
  \GG_N=\bigl\{(A^*)^n g_j:\ j\in J,\ 0\le n\le N\bigr\},
\end{equation}
and stable recovery of every $f$ holds precisely when $\GG_N$ is a frame for $\Hi$. The aim is to trade spatial measurements for repeated measurements in time. Frame characterisations for operator iterations are developed in \cite{ACMT2017,CHP2020,ACNP2025,CMPP2020,CMS2023,ACCMP2017normal,AP2017}.
There is also a systems-theoretic interpretation \cite{ACKM2026survey}, $\Ob_N$ is the observability map of $x_{n+1}=Ax_n$, $y_n=Cx_n$ with $C=G^*$, and the frame inequality is finite-horizon exact observability. See \cite{ADK2013,ADK2015} for the origins of the space--time trade-off and \cite{DMM2021control} for the control-theoretic dictionary.

Model \eqref{eq:classical-data} presupposes a known evolution operator. In applications a nominal operator may describe only part of the evolution, with unknown additive inputs accounting for external sources or unresolved effects. Such inputs accumulate over time and contaminate the orbit samples used to reconstruct the initial state. We retain the nominal operator $A$ and study this additive-forcing problem. An uncertain transition operator can also be written algebraically in this form, but the resulting state-dependent forcing need not satisfy the Gaussian assumptions used below.

Two branches of the literature respond to this situation, neither of them quite in the way we do.  The inhomogeneous model $x(n+1)=Ex(n)+w(n)$ has been studied extensively, but as \emph{source recovery}, there $w$ is a physical forcing, like pollutant emissions, diffusive or burst-like sources, which is itself the object of interest and is usually given structure, such as taking values in a known subspace or being periodic in time \cite{ADGMM2023,ACM2025periodic,ADMM2024,AHKK2023burst}. Uncertainty in the evolution operator, by contrast, is addressed by \emph{system identification}, which estimates the spectrum of $E$ from the same space--time samples by Krylov and Prony-type methods \cite{AK2016krylov,Tang2017sysid}, or learns the transition operator directly \cite{KMT2024learning}. Perturbation and error analyses of the initial-state problem are carried out in \cite{DMO2023error}.  Our objective is recovery of the initial state in the presence of a  forcing. We do not identify the evolution operator. Instead, we estimate the forcing under a specified covariance model and use it to correct the samples before reconstructing the initial state.

An alternative way is to think the problem as an unknown input, with the state as the object of interest.  Written as $x(n+1)=Ax(n)+w(n)$, $y_n=Cx(n)$ with $C=G^*$, this is \emph{state estimation with unknown input}, developed through unknown-input observers and unbiased minimum-variance filters \cite{Kitanidis1987,DarouachZasadzinskiXu1994,DarouachZasadzinski1997,HouMuller1992,HouPatton1998,Hsieh2000,GillijnsDeMoor2007}.  That theory treats
the input as a deterministic signal and pays for it with a rank condition on the pair (input, output). For an input entering through $G_w$ one requires $\operatorname{rank}(CG_w)=\operatorname{rank}(G_w)$, so that the input leaves a signature in the instantaneous measurement.  When the forcing can act anywhere in the state space, $G_w=I$, the condition reads $\operatorname{rank}G^*=\dim\Hi$, that is, the sensors would already have to span $\Hi$.  Dynamical sampling is the regime in which they do not.

Two responses to a deficient rank condition are available, and we take the
second.  One is to enlarge the output until the condition is restored, by
adjoining delayed or successively propagated measurements; this is a standard
device \cite{SundaramHadjicostis2007}, governed by the strong-observability
notions of \cite{Hautus1983}.  In that language the augmented system
\eqref{eq:double} and the depth-$M$ hierarchy of
Proposition \ref{prop:hierarchy} are output augmentation, and
\eqref{eq:residual-intro} is the signal-free residual it produces.  The other response, and the one developed here, is not to assume that the enlargement succeeds.  The forcing is given a prior, the residuals determine its component in the observed subspace, and the component that remains invisible is priced by a conditional covariance rather than declared unrecoverable (Section \ref{sec:lift-sub}).  The deterministic theory has no counterpart in that regime, because the rank condition is what it assumes.

We therefore let the observed trajectory be
\begin{equation}\label{eq:model}
  y_n=A^nf+\delta_n,\qquad \delta_0=0,
\end{equation}
with $\delta=(\delta_n)$ an unknown deviation, and we sample $y_n$ rather than $A^nf$.  Equivalently, writing $w_n:=\delta_n-A\delta_{n-1}$,
\begin{equation}\label{eq:recursion}
  y_n=Ay_{n-1}+w_n ,
\end{equation}
so that $w$ is the \emph{unmodeled forcing} and $\delta$ its accumulation along the nominal dynamics.

Without a normalisation at time zero, an unconstrained deviation can absorb any change of initial state (Section \ref{sec:ident}). The condition $\delta_0=0$ fixes this trajectory-level ambiguity. It does not, however, make the forcing observable from the original samples, which mix the initial-state orbit and the accumulated deviation.

The starting point of this work is to sample against the \emph{augmented system}
\begin{equation}\label{eq:double}
  \GG_1=\{g_j\}_{j\in J}\ \cup\ \{A^*g_j\}_{j\in J} .
\end{equation}
Indeed, \eqref{eq:recursion} gives
\begin{equation}\label{eq:residual-intro}
  \ip{y_n}{g_j}-\ip{y_{n-1}}{A^*g_j}
   =\ip{y_n-Ay_{n-1}}{g_j}=\ip{w_n}{g_j}.
\end{equation}
Whenever the auxiliary measurements are available, this identity cancels the initial-state term exactly. It exposes samples of the forcing and permits a two-stage procedure, estimate the forcing from the residuals, then correct the orbit samples and invert the frame. Neither full forcing recovery nor exact initial-state recovery follows from augmentation alone. Both still depend on the sampling geometry and, for the statistical estimator, the forcing prior. The nominal sample budget doubles when the two sensor banks are distinct.

With exact, signal-free residuals in hand we model the unmodeled forcing as a centered Gaussian vector over the observation window and estimate it from the residual data. The reconstruction then has a two-stage form, subtract the estimated deviation from the raw samples, and invert the analysis operator of the truncated dynamical sampling system. Each stage contributes one factor to the error, and the frame inequality is what converts a bound on the learning error into a bound on the reconstruction error.

Two features matter for the analysis. First, the residual likelihood contains no initial-state term under the additive-forcing model. Identifiable temporal hyperparameters can therefore be fitted without estimating $f$ first, although this does not identify an arbitrary spatial covariance from partial observations. Second, errors in estimating $w_k$ are propagated by $A^{n-k}$ and accumulate along the horizon, so the benefit of increasing the horizon must be assessed together with the frame geometry.

Model \eqref{eq:model} is, after the substitution \eqref{eq:recursion}, the inhomogeneous discrete-time model $x(n+1)=Ex(n)+w(n)$, $y_{n,j}=\ip{x(n)}{g_j}$ that organises the \emph{source recovery} branch of dynamical sampling \cite{ACKM2026survey}. That literature asks, under structural hypotheses on the source space $W\ni w(n)$, for $w$ constant in time, necessary and sufficient conditions on the sensors are given in \cite{ADMM2024} in finite dimensions and in \cite{ADGMM2023} in general, for $N$-periodic sources, stable recovery holds if and only if $\{P_W(T_s^{*}g_j)\}_{j}$ is a frame for $W$ for each $s$, where $T_s=(e^{2\pi is/N}I-A)^{-1}$ and $\opnorm A<1$ \cite{ACM2025periodic}. Burst-like continuous-time sources are treated in \cite{AHKK2023burst}. A recurring theme there is that $W=\Hi$ makes the problem easy, recovery holds iff $\{g_j\}$ is already a frame, and finitely many time samples suffice, while a proper subspace $W\subsetneq\Hi$ may force infinite time sampling \cite{ADGMM2023}.  Our visibility condition \eqref{eq:visibility} is the counterpart of that dichotomy in a probabilistic formulation, reached by a different route.

The distinguishing features are the initial-state reconstruction objective, the Gaussian covariance model for the forcing signal, and the enlarged sensor system that makes the residuals signal-free. The spatial covariance is assumed specified in the theory and determines the supported forcing directions. Fitting temporal parameters is a separate numerical extension. Theorem \ref{thm:main} gives an exact average-case reconstruction identity for fixed parameters, while Remark \ref{rem:relaxations} explains the losses introduced by coarser bounds.

The central question is how to recover the initial state when unknown
additive forcing contaminates the dynamical samples.  Our approach uses an
augmented measurement system to separate the two contributions before
reconstruction. Measurements against the original sensors and their adjoint
images yield one-step residuals whose observation equation contains the
forcing but no initial-state term. At greater augmentation depths, the same
construction provides observations against successive adjoint images of the
sensors. This gives a two-stage procedure: estimate the sampled accumulated
forcing from the residuals, then reconstruct the initial state from the
corrected orbit samples.

The residual protocol is output augmentation, and as such it is known.  The
contribution begins where augmentation stops working.  When the observed
subspace $S_{M-1}$ does not contain the directions in which the forcing acts,
the rank condition of the unknown-input literature fails and its estimators
are undefined; the question we answer is what remains recoverable, and at what
cost.  Under the separable prior we characterise the forcing uncertainty that
survives exact observation of every visible component, as a conditional
covariance with two-sided bounds through covariance-weighted subspaces
(Lemma \ref{lem:lift}).  Its trace $\kappa_{\mathrm{eff},M}$ is a continuous
relaxation of that rank condition, equal to zero exactly when the condition
holds (Remark \ref{rem:rank}) and positive, with a computable value, when it
does not.  Two consequences are worth separating: forcing supported outside
the observed subspace need not be irreducible, since prior correlations may
make it predictable from the visible components; and a positive
forcing-estimation floor does not automatically imply a positive
reconstruction floor, since the sampling and reconstruction operators can
annihilate some error directions.

Gaussian conditioning and canonical dual-frame reconstruction
are established tools. The contribution here is their
integration with this residual measurement protocol, together
with an analysis that retains the sampling geometry and the
dependence induced by shared measurement noise. In particular,
the forcing-estimation error and the noise in the samples
used for reconstruction are generally correlated, even when
the raw sensor-bank noises are independent. For a correctly
specified Gaussian model with fixed covariance parameters,
we derive an exact unconditional reconstruction covariance
and second-moment identity that include this dependence.
The resulting expression can be evaluated from the nominal
operator, sensors, and covariance model before observing
the initial state.

The analysis keeps two geometric requirements apart, and they are easily
conflated.  Augmentation depth determines which forcing directions the
residuals observe; the base-bank orbit determines the stability of nominal
initial-state reconstruction, through the frame bound $\alpha_N$.  The first
is a statement about $S_{M-1}$ and the prior, the second about $\GG_N$ and the
dynamics, and neither implies the other.  The reconstruction identity is where
they meet, since it carries both the conditional covariance of the forcing
error and the dual-frame operator that propagates it.

These results concern a specified spatial covariance and fixed statistical
parameters. They do not establish recovery of an unknown forcing subspace, nor
optimality among all estimators using the augmented data.  Fitting the
covariance parameters from the same residuals is a data-dependent procedure
that lies outside the fixed-parameter Gaussian guarantees; its status is
discussed in Remark \ref{rem:empirical-bayes}.

The work is organized as follows. Sections \ref{sec:double} and \ref{sec:ident} introduce the
residual protocol and distinguish trajectory identifiability
from nominal recoverability. Section \ref{sec:gp} develops
the Gaussian estimator and conditional uncertainty analysis.
Sections \ref{sec:guarantee} derives the
reconstruction identity, its relaxations, and their horizon
dependence. The experiments in Section \ref{sec:numerics}
validate the covariance calculation and illustrate the
different roles of depth, horizon, and covariance structure.

\section{Preliminaries}\label{sec:setting}

Let $J$ be a countable set and $\{g_j\}_{j\in J}$ a Bessel sequence with bound $\beta>0$, i.e.\ $\sum_{j}|\ip{u}{g_j}|^2\le\beta\norm{u}^2$ for all $u\in\Hi$.  We write $|J|=p$ when $J$ is finite.

We write $G^*:\Hi\to\ell^2(J)$, $G^*u=\bigl(\ip{u}{g_j}\bigr)_{j\in J}$, for the analysis operator of $\{g_j\}$, and $S:=\overline{\operatorname{span}}\, \{g_j\}_{j\in J}$, so that $G^*$ is injective on $S$.

\begin{definition}\label{def:frame}
For $N\in\N$ let $\GG_N$ be as in \eqref{eq:frame-system} and let $\Ob_N:\Hi\to\ell^2(J\times\{0,\dots,N\})$, $(\Ob_Nu)_{j,n}=\ip{u}{(A^*)^ng_j}$, be its analysis operator.  We say $\GG_N$ is a frame with bounds $0<\alpha_N\le\beta_N$ if $\alpha_N\norm{u}^2\le\norm{\Ob_Nu}^2\le\beta_N\norm{u}^2$ for all $u\in\Hi$.  We write $\Syn_N$ for the synthesis operator of the canonical dual frame $\{\tilde g_{j,n}\}$, so that $\Syn_N\Ob_N=I$ and $\opnorm{\Syn_N}\le\alpha_N^{-1/2}$.
\end{definition}

In finite dimensions, with $C=G^*$, the observability matrix is
\[
  \Ob_N=\begin{bmatrix}C\\CA\\\vdots\\CA^N\end{bmatrix}.
\]
If this matrix has full column rank, its optimal lower frame bound is $\alpha_N=s_{\min}(\Ob_N)^2=\lambda_{\min}(\Ob_N^*\Ob_N)$. Thus the numerical lower frame bound is a squared singular value.

We write $P:=P_S$ for the orthogonal projection onto $S$ and $Q:=I-P$. We distinguish two subspace angles used below.

\begin{definition}\label{def:angle}
For closed subspaces $\mathcal U,\mathcal V$ of a Hilbert space, the
\emph{Dixmier cosine} is
\[
  c_D(\mathcal U,\mathcal V)=\sup\bigl\{|\ip{u}{v}|:
    u\in\mathcal U,\ v\in\mathcal V,\ \norm u,\norm v\le1\bigr\}.
\]
The \emph{Friedrichs cosine} removes the intersection:
\[
  c(\mathcal U,\mathcal V)=c_D\bigl(
    \mathcal U\ominus(\mathcal U\cap\mathcal V),
    \mathcal V\ominus(\mathcal U\cap\mathcal V)\bigr).
\]
These agree when $\mathcal U\cap\mathcal V=\{0\}$, otherwise $c_D=1$
whenever the intersection is nonzero. The condition $c<1$ is equivalent to closedness of $\mathcal U+\mathcal V$, see \cite{Deutsch1995}.
\end{definition}

\begin{assumption}\label{ass:realisable}
The sensors $A^*g_j$, $j\in J$, are physically realisable, i.e.\ the quantities $\ip{y}{A^*g_j}=\ip{Ay}{g_j}$ can be measured.
\end{assumption}

Assumption \ref{ass:realisable} is a restriction and must be checked case by case, though it is benign in the standard families. For instance, if $A$ is the shift by $s$ and $g_j$ the evaluation at $jm$, then $A^{*}g_j$ is the evaluation at $jm-s$, an identical sensor placed upstream.
When $m\mid s$ one has $A^{*}g_j=g_{j-s/m}$, the auxiliary sensors are already present, and an array deployed along the flow direction supplies the depth-$M$ samples for free, the upstream sensor is the auxiliary sensor of its downstream neighbour.  When $m\nmid s$ the auxiliary sites interleave the existing ones and must be added; they remain of the same physical type, but the budget $p(M+1)(N+1)$ of Proposition \ref{prop:hierarchy} is then paid in full.

Sensor reuse does not imply nominal recoverability. For a pure integer shift on $\ell^2(\Z)$ with point sensors on $m\Z$, the iterates visit all residue classes modulo $m$ precisely when $\gcd(s,m)=1$; in that case $N\ge m-1$ suffices. If $m\mid s$ and $m>1$, every iterate remains on the original sublattice, so no horizon gives a frame for the full space. A filtered transport model requires a separate observability check.

If $A$ is a convolution by a kernel $a$, so that $Ax=a*x$, and $g_j$ is a mollified point evaluation, then $A^{*}g_j=\tilde a*g_j$ with $\tilde a(t)=\overline{a(-t)}$, the same evaluation seen through a further filter, that is, a sensor at the same location with a wider footprint.  For self-adjoint $A$ this is simply $Ag_j$.  Throughout, $\hat a$ denotes the Fourier transform of the kernel, the \emph{symbol} of $A$, so that $A$ acts as multiplication by $\hat a$ on the Fourier side; $A$ is self-adjoint exactly when $\hat a$ is real-valued, and $\opnorm A=\sup_\xi|\hat a(\xi)|$.

\section{Augmented sampling and exact residuals}\label{sec:double}

At each time $n$ we record the two sample vectors
\begin{equation}\label{eq:samples}
  u_n:=G^*y_n\in\ell^2(J),\qquad
  v_n:=G^*Ay_n=\bigl(\ip{y_n}{A^*g_j}\bigr)_j\in\ell^2(J).
\end{equation}

Let $y$ satisfy \eqref{eq:recursion}.  Then, for every $n\ge1$, the residual
\begin{equation}\label{eq:residual}
  z_n:=u_n-v_{n-1}=G^*w_n
\end{equation}
does not depend on $f$.  Indeed, $u_n-v_{n-1}=G^*y_n-G^*Ay_{n-1}=G^*(y_n-Ay_{n-1})=G^*w_n$ by \eqref{eq:recursion}, and $w_n$ is a function of $\delta$ alone.

The elimination of $f$ is what makes the whole scheme work. In the raw data $u_n$, signal and deviation are mixed, whereas the residual $z_n$ sees only the deviation. Learning may therefore be carried out first, on data whose observation equation contains no initial-state term, and reconstruction afterwards.

The scheme doubles the sample budget, from $p(N+1)$ to $2p(N+1)$.  It does not, by itself, resolve the identifiability problem (treated in Section \ref{sec:ident}), the residuals are invariant under the ambiguity $(f,\delta_n)\mapsto(f',\delta_n+A^n(f-f'))$, because that substitution leaves $w_n$ unchanged. Exact residuals give access to the forcing; they do not tell us how to split the data between signal and deviation.

The above is the case $M=1$ of a hierarchy, the following Proposition extends this idea.

\begin{proposition}[Depth-$M$ sampling]\label{prop:hierarchy}
Fix $M\ge1$ and suppose that at each time $n$ one records the $M+1$ sample vectors
\[
  u^{(m)}_n:=\bigl(\ip{y_n}{(A^*)^mg_j}\bigr)_{j\in J},
  \qquad 0\le m\le M ,
\]
so that $u^{(0)}_n=u_n$ and $u^{(1)}_n=v_n$ in the notation of \eqref{eq:samples}.  Then the following hold.
\begin{enumerate}
  \item For every $n\ge1$ and every $0\le m\le M-1$,
        \begin{equation}\label{eq:hier-residual}
          u^{(m)}_n-u^{(m+1)}_{n-1}
            =\bigl(\ip{w_n}{(A^*)^mg_j}\bigr)_{j\in J} ,
        \end{equation}
        so that the first $M$ moments of the forcing against the sensors are measured exactly, and each of them is free of $f$.
  \item Consequently, for every $n\ge1$ and every $j\in J$,
        \begin{equation}\label{eq:delta-samples}
          \ip{\delta_n}{g_j}
          =\underbrace{\sum_{m=0}^{\min(M,n)-1}
              \bigl(u^{(m)}_{n-m}-u^{(m+1)}_{n-m-1}\bigr)_j}
              _{\text{measured exactly}}
          \;+\;\underbrace{\sum_{m=M}^{n-1}
              \ip{w_{n-m}}{(A^*)^{m}g_j}}_{\text{must be estimated}} ,
        \end{equation}
        the second sum being empty when $n\le M$.
\end{enumerate}
The sample budget over the horizon $0\le n\le N$ is $|J|(M+1)(N+1)$.
\end{proposition}

\begin{proof}
(1) Both vectors on the left of \eqref{eq:hier-residual} are recorded, since $m\le M-1$ forces $m+1\le M$.  By definition of the adjoint,
\[
  \ip{y_{n-1}}{(A^*)^{m+1}g_j}
  =\ip{y_{n-1}}{A^*\bigl((A^*)^{m}g_j\bigr)}
  =\ip{Ay_{n-1}}{(A^*)^{m}g_j},
\]
so that, subtracting from $u^{(m)}_n=\ip{y_n}{(A^*)^mg_j}$ and using \eqref{eq:recursion},
\[
  \bigl(u^{(m)}_n-u^{(m+1)}_{n-1}\bigr)_j
  =\ip{y_n-Ay_{n-1}}{(A^*)^{m}g_j}
  =\ip{w_n}{(A^*)^{m}g_j}.
\]
The right-hand side is a function of $\delta$ alone, so the signal has cancelled; note that no property of $A$, of the sensors or of $\delta$ has been used beyond the availability of the measurements themselves.

(2) Unrolling \eqref{eq:recursion} from $\delta_0=0$ gives $\delta_n=\sum_{k=1}^{n}A^{\,n-k}w_k$, whence, moving $A$ across the inner
product,
\[
  \ip{\delta_n}{g_j}=\sum_{k=1}^{n}\ip{A^{\,n-k}w_k}{g_j}
   =\sum_{k=1}^{n}\ip{w_k}{(A^*)^{\,n-k}g_j}.
\]
Reindex by the order $m:=n-k$ of the moment, so that $k=n-m$ and $m$ runs over $0,\dots,n-1$:
\[
  \ip{\delta_n}{g_j}=\sum_{m=0}^{n-1}\ip{w_{n-m}}{(A^*)^{m}g_j}.
\]
The term of index $m$ is exactly the right-hand side of \eqref{eq:hier-residual} with $n$ replaced by $n-m$.  It is therefore measured precisely when that instance of \eqref{eq:hier-residual} is available, namely when $m\le M-1$ (the moment is recorded) and $n-m\ge1$ (the time index is admissible); the two conditions read $m\le\min(M,n)-1$.  For such $m$,
\[
  \ip{w_{n-m}}{(A^*)^{m}g_j}
   =\bigl(u^{(m)}_{n-m}-u^{(m+1)}_{n-m-1}\bigr)_j ,
\]
and splitting the sum at $m=\min(M,n)-1$ gives \eqref{eq:delta-samples}.  The remaining terms have $m\ge M$, they involve moments of order at least $M$, which are not recorded.  When $n\le M$ every term is of the first kind and the second sum is empty, so the whole of $\ip{\delta_n}{g_j}$ is measured.

Finally, at each of the $N+1$ times one records $M+1$ vectors with $|J|$ entries each, which is the stated budget.
\end{proof}

Depth-$M$ sampling trades sensor budget against the portion of the deviation that must be inferred rather than measured.  With depth $M$ the residuals determine the components of $w_n$ in
\begin{equation}\label{eq:SM}
  S_{M-1}:=\overline{\operatorname{span}}\bigl\{(A^*)^mg_j:\ j\in J,\ 0\le m\le M-1\bigr\},
\end{equation}
in the noiseless case. Inference about its orthogonal complement requires additional prior information.  This subspace grows with the depth, $\dim S_{M-1}\le M|J|$, so the depth is exactly what gives visibility of the forcing; the horizon $N$, by contrast, supplies more residuals, but each of them equally partial.  The observation suggests a design criterion, the sampling system should be chosen so as to \emph{observe the model error}, not only the signal.

Both ends of the hierarchy are to be avoided, for opposite reasons. At $M=N$ everything is measured exactly.  The time-$0$ measurements are then $\ip{f}{(A^*)^mg_j}$ for $m\le N$, that is, the full analysis $\Ob_Nf$, so if $\GG_N$ is a frame the signal is recovered from a single instant and the dynamics plays no role.  Since $\delta_0=0$ those measurements are moreover uncontaminated, this end of the hierarchy does not merely dispense with the evolution, it sidesteps the problem of this paper entirely.  It does so by demanding $N+1$ banks of sensors instead of one, which is precisely the hardware that dynamical sampling exists to avoid, and by spending $p(N+1)^2$ samples where $p(N+1)$ would do.

At $M=1$ the residuals give only $\ip{w_n}{g_j}$, the two sensor banks yield a single moment of the forcing, the second being consumed in forming the residual, so that only the component of $w_n$ in $S_0=S$ is determined.  Observable temporal parameters may still be estimated, but this does not by itself suffice for correcting the samples, since \eqref{eq:delta-samples} also requires $\ip{w_k}{(A^*)^{m}g_j}=\ip{A^{m}w_k}{g_j}$ for $m\ge1$, and $A$ need not preserve $S^{\perp}$, an invisible forcing component can reappear in the observed directions at later times.  How much this costs depends on where the forcing acts, and is quantified in Section \ref{sec:lift-sub}.

The useful regime is therefore intermediate, $M$ small compared with $N$, so that the dynamics still does the work, but large enough that $S_{M-1}$ contains the directions along which the forcing actually acts.

\section{Identifiability}\label{sec:ident}

Suppose the normalisation $\delta_0=0$ in \eqref{eq:model} is dropped, so that $\delta$ ranges over all sequences in $\Hi$.  Then for all $f,f'\in\Hi$ and every $\delta$, the pairs $(f,\delta)$ and $(f',\delta')$ with $\delta'_n:=\delta_n+A^n(f-f')$ produce the same trajectory, hence the same data \eqref{eq:samples} at every depth $M$.  No amount of data distinguishes $f$ from $f'$. In fact, $A^nf'+\delta'_n=A^nf'+\delta_n+A^nf-A^nf'=A^nf+\delta_n=y_n$, so the trajectories, hence all samples, coincide.

\begin{proposition}\label{cor:identifiable}
Under $\delta_0=0$, if $(f,\delta)$ and $(f',\delta')$ both satisfy \eqref{eq:model} with $\delta_0=\delta'_0=0$ and produce the same trajectory $y$, then $f=f'$.
\end{proposition}
\begin{proof}
Both pairs produce the same trajectory, so in particular $y_0$ is the same for both.  Evaluating \eqref{eq:model} at $n=0$ for the first pair gives
\[
  y_0=A^0f+\delta_0=f ,
\]
since $A^0=I$ and $\delta_0=0$, and the same computation for the second pair gives $y_0=f'$.  Hence $f=y_0=f'$.
\end{proof}

Proposition \ref{cor:identifiable} says that $f$ is determined by the \emph{trajectory}, it does not say that $f$ can be computed from the \emph{samples}.  What the time-$0$ data deliver is $u_0=G^*y_0=G^*f$, which determines only $P_Sf$, and in the regime of interest (few sensors, so $S\subsetneq\Hi$) that is strictly less than $f$.

Stable inversion of the \emph{nominal orbit samples} is a separate property, it holds when $\GG_N$ is a frame, that is, when $\alpha_N>0$ in Definition \ref{def:frame}. We refer to this as nominal recoverability. It does not ensure exact recovery from samples contaminated by arbitrary forcing. Even with $\delta_0=0$, unobserved forcing components can conceal a change of initial state in the partial measurements.

Concretely, the raw samples are
\[
  u_n=G^*y_n=(\Ob_Nf)_{\cdot,n}+G^*\delta_n.
\]
Reconstruction therefore requires both a correction for the sampled deviation and stable inversion of $\Ob_N$. The frame condition provides the second; the residuals and the prior control the error in the first. For few sensors, the static system need not span $\Hi$, while its iterates can form a frame.

\begin{remark}\label{rem:friedrichs-static}
There is a purely deterministic way of measuring what the unknown deviation costs. Let $\HH\subseteq\ell^2(\{0,\dots,N\};\Hi)$ be a closed class of admissible deviations, let $(\Sa_N\delta)_{j,n}=\ip{\delta_n}{g_j}$, and set $\mathcal U:=\overline{\ran}\,\Ob_N$, $\mathcal V:=\overline{\ran}\,(\Sa_N|_{\HH})$.  If $\mathcal U\cap\mathcal V=\{0\}$ and $c:=c(\mathcal U,\mathcal V)<1$, then
\begin{equation}\label{eq:oblique-lower}
  \norm{\Ob_Nf+\Sa_N\delta}^2\ \ge\ (1-c^2)\norm{\Ob_Nf}^2
   \ \ge\ (1-c^2)\,\alpha_N\norm{f}^2
   \qquad(f\in\Hi,\ \delta\in\HH),
\end{equation}
with sharp constant, the reconstruction is the oblique projection onto $\mathcal U$ along $\mathcal V$, of norm $(1-c^2)^{-1/2}$, in the sense of consistent sampling \cite{Eldar2003}.  The effective stability constant is $\alpha_N(1-c^2)$, so the cost of the unknown deviation is exactly the alignment between what a signal and what an admissible deviation can produce in the data.  We do not build on \eqref{eq:oblique-lower}, because it assumes a class $\HH$ with computable $c<1$, whereas an unrestricted $\HH$ gives $\mathcal U\subseteq\mathcal V$, violating the required trivial intersection by the ambiguity at the start of this section; and it makes no use of the residuals \eqref{eq:residual}, which say far more than that $\delta$ lies in some subspace.
\end{remark}

\section{The Gaussian model}\label{sec:gp}

Proposition \ref{cor:identifiable} settles identifiability but leaves the estimation problem empty.  If $(w_n)$ is an arbitrary sequence in $\Hi$, each residual $z_n$ is one partial observation of one unknown vector $w_n$ and constrains no other, the horizon $N$ brings $N$ new data and $N$ new unknowns, and the component of $w_n$ outside $S_{M-1}$ is never seen at all.  Something must couple the $w_n$ to one another, and we let that something be a prior covariance.

\begin{assumption}\label{ass:prior}
The forcing over the observation window is centered and jointly Gaussian,
\[
  w=(w_1,\dots,w_N)\sim\mathcal N(0,\mathcal K),
\]
with $\mathcal K$ trace class on $\bigoplus_{n=1}^{N}\Hi$. The initial state $f$ is fixed; the measurement noises are jointly centered Gaussian and independent of $w$, with trace-class covariance in their observation spaces. In a complex Hilbert space the Gaussian vectors are assumed proper (their pseudo-covariances vanish. 
Alternatively all statements can be read in a real Hilbert space. Covariance and noise parameters are fixed, independently of the residuals used for reconstruction. We use the separable choice
\begin{equation}\label{eq:sep-prior}
  \Cov(w_n,w_m)=K(n,m)\,C_w ,
\end{equation}
where $K$ is a positive-definite kernel on $\{1,\dots,N\}$, normalised so that $K(n,n)=1$, and $C_w$ is self-adjoint, positive and trace class on $\Hi$.
\end{assumption}

Trace-class $C_w$ is what makes the forcing take values in $\Hi$ almost surely \cite{Stuart2010}, and \eqref{eq:sep-prior} separates the two questions that the rest of the paper keeps apart, $K$ says how the forcing is correlated in time, which is what makes one residual informative about another, and $C_w$ says where in $\Hi$ it acts, which is what decides how much of it the sensors see.  
The construction is the Hilbert space analogue of a linear latent force model \cite{ALL2013}, see \cite{RW2006} for the underlying machinery and \cite{AHKLLV2018noise} for additive measurement noise in dynamical sampling.

The prior correlates the forcing across the window, which is what makes a residual at one time informative about the forcing at another. 

By \eqref{eq:residual}, the quantities
\begin{equation}\label{eq:dataset}
  \D_N=\{z_1,\dots,z_N\},\qquad
  z_n=u_n-v_{n-1}=G^*w_n+\varepsilon^u_n-\varepsilon^v_{n-1},
\end{equation}
are available from the measurements (here $u_n,v_n$ denote their noisy
versions), and each $z_n$ is a noisy bounded linear functional of the unknown Gaussian vector $w$ of Assumption \ref{ass:prior}, and the estimation problem is Gaussian conditioning rather than regression.  Two properties of $\D_N$ matter;  its observation equation is free of $f$, which permits fitting identifiable
temporal parameters from residuals alone; and it is causal, since $z_n$ uses only the times $n$ and $n-1$, so it may be accumulated online, offline from earlier runs, or both.

\begin{remark}\label{rem:noise}
For finite $J$, if the two sensor banks are corrupted by noises independent across banks and times, $\varepsilon^u_n,\varepsilon^v_n\sim\mathcal N(0,\sigma^2I)$, then $z_n$ has noise covariance $2\sigma^2I$.  Moreover $z_n$ and $z_{n'}$ involve disjoint noise terms for $n\ne n'$, since the first uses $(\varepsilon^u_n,\varepsilon^v_{n-1})$ and the second $(\varepsilon^u_{n'},\varepsilon^v_{n'-1})$; the noise covariance of the whole residual vector is therefore exactly $2\sigma^2I$, with no temporal structure.  This is a specific benefit of measuring the two sample banks independently rather than propagating one from the other. For countably infinite $J$, $\mathcal N(0,\sigma^2I)$ is not an $\ell^2(J)$-valued random vector. The Hilbert-valued theory instead requires trace-class observation-noise covariance. 

At depth $M>1$, independence of the raw sensor-bank noises does not imply independence of all residual noises. If $e_n^{(m)}$ is the raw noise in bank $m$, the residual noise is $e_n^{(m)}-e_{n-1}^{(m+1)}$, adjacent residuals can share a raw noise term. The full covariance must therefore be retained.
\end{remark}

\subsection{The posterior}

The residuals give $G^*w_n$, not $w_n$, whereas \eqref{eq:delta-samples} at depth $M=1$ also needs moments of orders $m\ge1$. For general depth $M$, let $\Sa$ act at each time as the analysis operator of $\GG_{M-1}$, and stack the residuals as $z=\Sa w+\eta$. At depth $M$, $\D_N$ denotes all residuals in \eqref{eq:hier-residual} over the window. Write $\Sigma_\eta=\Cov(\eta)$ and $D_z=\Sa\mathcal K\Sa^*+\Sigma_\eta$. For finite $J$, suppose $D_z$ is positive definite; this holds if $\Sigma_\eta$ is positive definite. Then Gaussian conditioning gives
\begin{align}
  \hat w&=\E[w\mid z]=Fz,\qquad
  F:=\mathcal K\Sa^*D_z^{-1}, \label{eq:post-mean}\\
  \Sigma_r&:=\Cov(w-\hat w)
    =\mathcal K-\mathcal K\Sa^*D_z^{-1}\Sa\mathcal K.
    \label{eq:Sigma-r}
\end{align}
At depth $M=1$ with the independent isotropic noise of Remark \ref{rem:noise}, $\Sigma_\eta=2\sigma^2I$. At higher depth the covariance includes the shared noise terms described there. In a noiseless finite observation model, $D_z^{-1}$ is replaced by its Moore--Penrose inverse.

For countably infinite $J$, trace-class noise does not in general make $D_z$ boundedly invertible. Conditional expectations and conditional covariances are still well defined, but the inverse formula must be read through Gaussian projection or a convergent finite-observation approximation. Whenever the displayed gain is used as an operator, we assume that its composition has a bounded extension. We do not assume an inverse exists on the whole infinite-dimensional observation space.

With fixed parameters, $\hat w$ is a linear Gaussian projection and $\Sigma_r$ does not depend on the realised residuals. It is both the conditional covariance of $w$ given $z$ and the unconditional covariance of $w-\hat w$, and is trace class.

\begin{remark}\label{rem:empirical-bayes}
If parameters are estimated from these same residuals by marginal likelihood, the gain becomes $F_{\hat\theta(z)}$. The resulting estimator is generally nonlinear in $z$, its error need not be Gaussian, and the plug-in posterior covariance is not its unconditional error covariance. The Gaussian guarantees below concern fixed parameters and do not apply directly to this empirical-Bayes procedure. Signal cancellation remains exact. The pathwise error identity also remains valid, but an unconditional second moment then requires the actual error law, including any nonzero bias. Parameters fitted from independent training data can instead be treated as fixed conditionally on those data; using the fitted covariance as the true covariance still requires a correctly specified model.
\end{remark}

\subsection{\texorpdfstring{The lift from $G^*w$ to $w$}{The lift from observed to full forcing}}\label{sec:lift-sub}

Noiseless residuals determine $Pw$ but not directly $Qw$, its projection onto $S^\perp$. Estimating the full forcing is one way to obtain the unobserved moments required in \eqref{eq:delta-samples}. The following isolates what the prior supplies. Throughout this subsection, $\{g_j\}$ is assumed to be a frame for $S$, so that exact knowledge of $G^*w$ is equivalent to that of $Pw$.

\begin{lemma}\label{lem:lift}
Let $w\sim\mathcal N(0,C_w)$ with $C_w$ positive and trace class, and suppose $Pw$ is known. Set $T:=C_w^{1/2}$ and let $\Pi$ be the orthogonal projection onto $\mathcal A:=\overline{\ran}(TP)$. The conditional law of $Qw$ is Gaussian with covariance
\begin{equation}\label{eq:schur-alt}
  \Sigma_\perp=QT(I-\Pi)TQ.
\end{equation}
In finite dimensions this is the Schur complement
\begin{equation}\label{eq:schur}
  \Sigma_\perp=QC_wQ-QC_wP(PC_wP)^\dagger PC_wQ.
\end{equation}
In infinite dimensions \eqref{eq:schur-alt} defines the covariance; the generalised-inverse expression is interpreted through this bounded
projection formula. Moreover,
\begin{equation}\label{eq:lift-twosided}
  (1-c_w^2)QC_wQ\preceq\Sigma_\perp\preceq QC_wQ,
\end{equation}
where
\begin{equation}\label{eq:cw}
  c_w:=c_D\bigl(\overline{\ran}(C_w^{1/2}P),
                     \overline{\ran}(C_w^{1/2}Q)\bigr)\in[0,1].
\end{equation}
This is a Dixmier cosine. It agrees with the Friedrichs cosine only when the two closed ranges have trivial intersection.
\end{lemma}

\begin{proof}
Gaussian conditioning is orthogonal projection onto the closed Gaussian inear span generated by $Pw$. Factoring the covariance as $T^2$, that span corresponds to $\overline{\ran}(TP)$ and its orthogonal complement gives
\eqref{eq:schur-alt}. In finite dimensions,
$TP(PC_wP)^\dagger PT=\Pi$, giving \eqref{eq:schur}.
The upper bound follows from $0\preceq I-\Pi\preceq I$. For the lower bound we work on quadratic forms, avoiding inverses of trace-class covariance operators.

Now fix $v\in\Hi$ and set $\eta:=TQv$, so that $\eta\in\mathcal B:=\overline{\ran}\,(TQ)$ and, by \eqref{eq:schur-alt},
\[
  \ip{\Sigma_\perp v}{v}=\norm{(I-\Pi)\eta}^2,
  \qquad \ip{QC_wQ\,v}{v}=\norm{\eta}^2 .
\]
By \eqref{eq:cw},
$\norm{\Pi\eta}=\sup\{|\ip{\eta}{u}|:u\in\mathcal A,\ \norm u\le1\}
 \le c_w\norm\eta$ for every $\eta\in\mathcal B$, whence
\[
  \ip{\Sigma_\perp v}{v}=\norm{\eta}^2-\norm{\Pi\eta}^2
   \ \ge\ (1-c_w^2)\norm{\eta}^2=(1-c_w^2)\ip{QC_wQ\,v}{v} .
\]
As $v\in\Hi$ was arbitrary, this is the left-hand inequality in \eqref{eq:lift-twosided}.

The upper bound is an equality when $C_w$ is block diagonal for $\Hi=S\oplus S^{\perp}$. If $QC_wQ\ne0$, the lower multiplier $1-c_w^2$ is sharp as a uniform relative quadratic-form bound, the range of $TQ$ is dense in $\mathcal B$, so normalised vectors in that range can approach the supremum defining $c_w$. This does not assert equality of the two covariance operators.
\end{proof}

Notice that $\Sigma_\perp=0$ precisely when $Qw$ is almost surely determined by $Pw$, the invisible part of the forcing is then a measurable function of the visible one and costs nothing.  The degenerate instance is $QC_wQ=0$, i.e.
\begin{equation}\label{eq:visibility}
  \ran C_w\subseteq S.
\end{equation}
This is the visibility condition, the forcing is supported in the sensor subspace and there is no invisible part. By the lower bound in
\eqref{eq:lift-twosided}, $c_w=1$ is necessary for $\Sigma_\perp$ to vanish without \eqref{eq:visibility}, but it is not sufficient, it expresses perfect or arbitrarily close alignment in some directions, not predictability of every invisible component.
Condition \eqref{eq:visibility} is a checkable design condition on the pair (prior, sensors), and $\kappa^2:=\tr(QC_wQ)$, the prior mass of the forcing in the directions invisible to the sensors, is its quantitative relaxation.

The discussion of \eqref{eq:SM} can now be made precise.  Under depth-$M$ sampling the residuals give the analysis of $w_n$ against $\GG_{M-1}$, hence determine $P_{S_{M-1}}w_n$ provided $\GG_{M-1}$ is a frame for $S_{M-1}$, and Lemma \ref{lem:lift} applies with $P$ replaced by $P_{S_{M-1}}$ and $Q$ by $Q_{M-1}:=I-P_{S_{M-1}}$.  The visibility condition \eqref{eq:visibility}
becomes $\ran C_w\subseteq S_{M-1}$, and the invisible prior mass is
\[
  \kappa_M^2=\tr\bigl(Q_{M-1}C_wQ_{M-1}\bigr),
\]
which is non-increasing in $M$ and vanishes exactly when $S_{M-1}\supseteq\ran C_w$.  For finite $J$, $\dim S_{M-1}\le M|J|$, so no depth below $\operatorname{rank}C_w/|J|$ can achieve this, the sensor budget of Proposition \ref{prop:hierarchy} and the effective dimension of the unmodeled dynamics are tied to one another.  Everything said so far for $P$ and $Q$ holds verbatim at depth $M$, with $c_{w,M}$ defined by \eqref{eq:cw} using the depth-$M$ projections. Let $\Sigma_{\perp,M}$ denote the corresponding conditional covariance and define
\begin{equation}\label{eq:kappa-eff}
  \kappa_{\mathrm{eff},M}^2:=\tr\Sigma_{\perp,M}.
\end{equation}
Then
$(1-c_{w,M}^2)\kappa_M^2\le\kappa_{\mathrm{eff},M}^2\le\kappa_M^2$. The effective mass is non-increasing with $M$, because conditioning on a
larger sensor subspace cannot increase covariance. It is positive exactly when there remains irreducible forcing uncertainty; visibility is sufficient, but not necessary, for it to vanish.

At depth $M=1$, a large invisible prior mass can limit forcing estimation unless prior correlations predict those components from the visible ones. Low-dimensional forcing supported near $S$ is one favourable situation, but its suitability is a modelling hypothesis rather than a consequence of augmentation alone.

\begin{remark}\label{rem:floor}
Under the separable prior, exact knowledge of every $Pw_m$ in the window leaves covariance $K\otimes\Sigma_\perp$ for $w$. To see this without unbounded generalised inverses, factor $\mathcal K$ as $(K^{1/2}\otimes T)^2$. Since $K$ is positive definite, the closed range associated with observing all visible components is $\mathbb C^N\otimes\overline{\ran}(TP)$ (or $\mathbb R^N$ in the real case). The unobserved projection therefore gives $K\otimes T(I-\Pi)T=K\otimes\Sigma_\perp$, where the last equality uses $(I-\Pi)TP=0$.

The residuals are functions of these visible components and independent measurement noise. They provide no more information than their exact values. Consequently, at depth $M$, the posterior error blocks obey
\[
  \tr\Cov(w_n-\hat w_n)\ge\tr\Sigma_{\perp,M}
     =\kappa_{\mathrm{eff},M}^2,
\]
because $K(n,n)=1$. The same lower bound holds for the mean-square forcing error of any estimator based on these residuals. It remains positive for all horizons when $\Sigma_{\perp,M}\ne0$, with a fixed spatial prior and depth. This is a forcing-estimation floor; it is not automatically a positive initial-state reconstruction floor, since the subsequent sampling and reconstruction operators can annihilate some error directions.
\end{remark}

\begin{example}\label{ex:predictability}Let $X,Y$ be independent standard real Gaussians. For $w=(X,X,Y)\in\R^3$ with $S=\operatorname{span}\{e_1\}$, observing $Pw$ determines the second coordinate but not the third. The covariance-weighted ranges in \eqref{eq:cw} intersect, so $c_w=1$, yet $\Sigma_\perp=\operatorname{diag}(0,0,1)\ne0$. In contrast, for $w=(X,X)\in\R^2$ with the same sensor direction, $\Sigma_\perp=0$ although $\ran C_w\not\subseteq S$. Thus neither failure of visibility nor $c_w=1$ alone determines the floor.\end{example}

\begin{remark}[The rank condition and its continuous relaxation]\label{rem:rank}
In finite dimensions $w$ lies almost surely in $\ran C_w$, so $Pw$ determines
$w$ exactly when $P$ is injective there, that is, when
\begin{equation}\label{eq:rank-condition}
  S^{\perp}\cap\ran C_w=\{0\}.
\end{equation}
Since $\Sigma_\perp=0$ precisely when $Qw$ is almost surely determined by
$Pw$, condition \eqref{eq:rank-condition} is equivalent to
$\kappa_{\mathrm{eff}}=0$.  It is also the rank condition
$\operatorname{rank}(CG_w)=\operatorname{rank}(G_w)$ of the unknown-input
literature, written for $C=G^*$ and the input factor $G_w=C_w^{1/2}$, whose
range is $\ran C_w$.  The visibility condition \eqref{eq:visibility} is
strictly stronger: $\ran C_w\subseteq S$ implies \eqref{eq:rank-condition},
and Example \ref{ex:predictability} separates the two, since $w=(X,X)$
satisfies \eqref{eq:rank-condition} while failing \eqref{eq:visibility}.

The Gaussian formulation does more than test \eqref{eq:rank-condition}: it
measures how far it fails.  The rank condition is binary, and the
deterministic estimators built on it are undefined once it is violated,
whereas $\Sigma_{\perp,M}$ is defined in either case and
Theorem \ref{thm:main} converts it into a reconstruction error.  In infinite
dimensions the argument above must be read through the Cameron--Martin space
of the prior rather than through $\ran C_w$, and \eqref{eq:schur-alt} is the
statement that survives.
\end{remark}

\subsection{A uniform bound on the forcing error}

The worst-case corollaries of Section \ref{sec:guarantee} need a single number controlling all the $r_k:=w_k-\hat w_k$ at once.  Since the window is finite, this costs only a union bound.

\begin{lemma}\label{lem:uniform}
Under the fixed-parameter Gaussian model, let $\Sigma_r$ be as in \eqref{eq:Sigma-r} and write, for its diagonal blocks, $\tau_k^2:=\tr\Cov(r_k)$ and $\varsigma_k^2:=\opnorm{\Cov(r_k)}$. Then for every $\delta\in(0,1)$,
\begin{equation}\label{eq:rho-total}
  \Prob\Bigl\{\ \norm{r_k}\le\tau_k
     +\varsigma_k\sqrt{2\log\bigl(N/(1-\delta)\bigr)}
     \ \text{ for all }1\le k\le N\ \Bigr\}\ \ge\ \delta ,
\end{equation}
and we write $\bar\rho_N$ for the largest of the $N$ bounds appearing there. By Remark \ref{rem:floor}, $\bar\rho_N\ge\max_k\tau_k \ge\kappa_{\mathrm{eff},M}$. Thus it has a positive, horizon-independent lower bound when $\Sigma_{\perp,M}\ne0$ and the spatial prior and depth are fixed. Failure of visibility alone does not imply this.
\end{lemma}

\begin{proof}
Each $r_k$ is a centered Gaussian vector in $\Hi$ with covariance $\Cov(r_k)$, and
\[
  \E\norm{r_k}\le(\E\norm{r_k}^2)^{1/2}=\tau_k.
\]
For $\varsigma_k>0$, Borell--TIS \cite{AdlerTaylor2007} applied to the $1$-Lipschitz map $v\mapsto\norm v$ gives $\Prob\{\norm{r_k}>\tau_k+t\}\le\exp(-t^2/2\varsigma_k^2)$.  Taking $t=\varsigma_k\sqrt{2\log(N/(1-\delta))}$ makes each of these probabilities at most $(1-\delta)/N$. If $\varsigma_k=0$, then $r_k=0$ almost surely. A union bound over $k=1,\dots,N$ gives \eqref{eq:rho-total}.
\end{proof}

The posterior covariance combines statistical uncertainty with the irreducible covariance of Remark \ref{rem:floor}. For nested windows from one consistent prior, additional observations can reduce the covariance at a fixed time, but neither a decrease of the maximum over an expanding window nor a rate of decrease is automatic. Invisible components can also benefit from observed correlations; it is their conditional remainder, rather than their entire prior variance, that cannot be removed.

\section{Reconstruction guarantee}\label{sec:guarantee}

The reconstruction proceeds in three steps, form the naive frame reconstruction from the raw samples, estimate the accumulated deviation from the estimated forcing, and subtract.

\begin{definition}\label{def:estimator}
Given $\D_N$ and the posterior mean $\hat w$, set $\hat\delta_n:=\sum_{k=1}^{n}A^{n-k}\hat w_k$ and
\begin{equation}\label{eq:estimator}
  \hat f:=\Syn_N\bigl(u-G^*\hat\delta\bigr),
  \qquad u=(u_n)_{n=0}^N,\quad G^*\hat\delta=(G^*\hat\delta_n)_{n=0}^N .
\end{equation}
\end{definition}

\subsection{An exact identity}

The error of \eqref{eq:estimator} is a \emph{linear} image of the forcing error and of the measurement noise, so its second moment can be written exactly, with no inequality at any step.  Let 
\begin{equation}\label{eq:Lacc}
  \Lacc:\ \bigoplus_{k=1}^{N}\Hi\ \longrightarrow\ \bigoplus_{n=1}^{N}\Hi,
  \qquad (\Lacc\,r)_n=\sum_{k=1}^{n}A^{\,n-k}r_k ,
\end{equation}
be the accumulation operator, let $\Sa_N$ act blockwise as $G^*$, and let $\iota_N$ embed the times $1,\dots,N$ into $0,\dots,N$.  Put
\begin{equation}\label{eq:TN}
  T_N:=\Syn_N\,\iota_N\,\Sa_N\,\Lacc,
  \qquad
  B_N:=\bigl[\,T_N\ \ \Syn_N\,\bigr].
\end{equation}

\begin{theorem}\label{thm:main}
Under Assumption \ref{ass:prior} with fixed parameters, assume $\GG_N$ is a frame for $\Hi$ with $\alpha_N>0$, let $\hat f$ be as in Definition \ref{def:estimator}, and write $r:=(w_k-\hat w_k)_{k=1}^N$ for the forcing error and $\varepsilon$ for the measurement noise on $u$.  Then
\begin{equation}\label{eq:exact-error}
  f-\hat f\ =\ -\,T_N\,r\ -\ \Syn_N\,\varepsilon .
\end{equation}
In particular the error does not depend on $f$, the estimator is unbiased. Moreover $(r,\varepsilon)$ is jointly Gaussian and centered, so with $\Sigma_{r\varepsilon}$ its covariance,
\begin{equation}\label{eq:exact-identity}
  \E\norm{f-\hat f}^2\ =\ \tr\bigl(B_N\,\Sigma_{r\varepsilon}\,B_N^{*}\bigr),
\end{equation}
the expectation being over the forcing and all measurement noises, with $f$ fixed.  Since $f-\hat f$ is then a centered Gaussian vector with covariance $\Sigma_f:=B_N\Sigma_{r\varepsilon}B_N^{*}$, for every $\delta\in(0,1)$
\begin{equation}\label{eq:exact-hp}
  \norm{f-\hat f}\ \le\ \sqrt{\tr\Sigma_f}
     +\sqrt{2\opnorm{\Sigma_f}\log(1/\delta)}
\end{equation}
with probability at least $1-\delta$.
\end{theorem}

\begin{proof}
From \eqref{eq:model}, $u=\Ob_Nf+\iota_N\Sa_N\delta+\varepsilon$ (using
$\delta_0=0$), while $\delta-\hat\delta=\Lacc\,r$ by \eqref{eq:Lacc}.  Hence
$u-G^*\hat\delta=\Ob_Nf+\iota_N\Sa_N\Lacc\,r+\varepsilon$, and applying
$\Syn_N$ with $\Syn_N\Ob_N=I$ gives \eqref{eq:exact-error}, in which $f$ does
not appear.  Under Assumption \ref{ass:prior} the model $z=\Sa w+\eta$ is linear and
Gaussian, so the posterior mean $\hat w$ is a Gaussian linear projection of $z$.
The residuals are linear in $w$ and the noises in all sensor banks, including
the auxiliary-bank noises; hence $r=w-\hat w$ is jointly Gaussian with
$\varepsilon$,
and it is centered because $\E[r]=\E\bigl[w-\E[w\mid\D_N]\bigr]=0$ by the tower
property.  Then \eqref{eq:exact-identity} is the trace formula for the second
moment of a linear image of a centered Gaussian vector, and \eqref{eq:exact-hp}
is Borell--TIS \cite{AdlerTaylor2007} applied to the $1$-Lipschitz map
$v\mapsto\norm v$, using $\E\norm{f-\hat f}\le(\tr\Sigma_f)^{1/2}$.  All
operators involved are bounded and $\Sigma_{r\varepsilon}$ is trace class,
hence so is $\Sigma_f$.
\end{proof}

\begin{remark}
\label{rem:joint}
The identity \eqref{eq:exact-identity} is unconditional. Given $z$, the
raw measurement noise generally has a nonzero conditional mean, so a
conditional mean-square formula must include the squared conditional bias.
The cross-covariance between $r$ and $\varepsilon$ also generally does not
vanish. With $F$ as in \eqref{eq:post-mean},
\begin{equation}\label{eq:joint-covariance}
  \Sigma_{r\varepsilon}=
  \begin{bmatrix}
    \Sigma_r & -F\Cov(\eta,\varepsilon)\\
    -\Cov(\varepsilon,\eta)F^* & \Sigma_\varepsilon
  \end{bmatrix},\qquad \Sigma_\varepsilon:=\Cov(\varepsilon).
\end{equation}
This follows from $r=(I-F\Sa)w-F\eta$ and independence of $w$ from the measurement noises. Omitting the off-diagonal blocks changes the predicted error and is not justified by independence of the raw sensor banks.
\end{remark}

\begin{remark}\label{rem:computable}
For fixed known covariance parameters, every ingredient of \eqref{eq:exact-identity} is available before measuring the signal, $T_N$ depends on $A$, the sensors and $N$, and $\Sigma_{r\varepsilon}$ follows from the forcing and raw noise covariances. The identity can therefore be used for sensor, depth and horizon design. For parameters fitted from the same data, the plug-in computation has a different status (Remark \ref{rem:empirical-bayes}). The abstract conditional expectation formulation allows countable $J$ with Hilbert-valued, trace-class noise; finite-matrix inverse formulas additionally require the conditions stated in the posterior subsection.
\end{remark}

\begin{remark}\label{rem:xi}
The reconstruction \eqref{eq:estimator} does not use $\hat w$ itself, only the sampled deviation $G^*\hat\delta$.  Writing
\[
  \xi:=\iota_N\Sa_N\Lacc\,w = G^*\delta
\]
for the sampled accumulated deviation, which is a linear image of $w$ and therefore Gaussian, the posterior mean commutes with that image, $\E[\xi\mid\D_N]=\iota_N\Sa_N\Lacc\,\hat w=G^*\hat\delta$, so the correction term of \eqref{eq:estimator} \emph{is} the posterior mean of $\xi$.  One may therefore condition $\xi$ on $z$ directly, without ever representing $\hat w$ in $\Hi$. The result is identical, and when $J$ is finite the computation takes place in $\ell^2(J)^{N+1}$ rather than in $\Hi^{N}$.

In that formulation \eqref{eq:exact-error} reads
$f-\hat f=-\Syn_N\zeta$ with $\zeta:=(\xi-\hat\xi)+\varepsilon$ the total
error in the sampled data, and \eqref{eq:exact-identity} becomes
\begin{equation}\label{eq:exact-xi}
  \E\norm{f-\hat f}^2=\tr\bigl(\Syn_N\Cov(\zeta)\,\Syn_N^{*}\bigr),
\end{equation}
which absorbs the cross-covariance of Remark \ref{rem:joint} rather than carrying it in a separate block.  We nevertheless keep the formulation in $w$, because it is the one in which Section \ref{sec:lift-sub} can be read, the floor of Remark \ref{rem:floor} is a statement about directions of $\Hi$ that the sensors never see, and in \eqref{eq:exact-xi} it survives only as a contribution to $\Cov(\zeta)$, with no geometric interpretation.
\end{remark}


When $\Sigma_{r\varepsilon}$ is not available, \eqref{eq:exact-identity} may be relaxed.  We record two levels of relaxation and quantify what each costs.

\begin{corollary}\label{cor:quadrature}
Under the hypotheses of Theorem \ref{thm:main}, suppose the centered errors $r_k$ are uncorrelated with $\tr\Cov(r_k)\le\rho^2$, and put
\begin{equation}\label{eq:Gammatilde}
  \tilde\Gamma_N:=\Bigl(\sum_{n=1}^{N}\sum_{m=0}^{n-1}
     \opnorm{A^{m}}^{2}\Bigr)^{1/2}.
\end{equation}
Then
\[
  \E\norm{f-\hat f}^2\ \le\ 2\alpha_N^{-1}
     \bigl(\beta\,\tilde\Gamma_N^{2}\,\rho^{2}+\E\norm{\varepsilon}^2\bigr),
\]
and $\E\norm{f-\hat f}^2\le\alpha_N^{-1}\beta\,\tilde\Gamma_N^{2}\rho^{2}$ in
the noiseless case $\varepsilon=0$.
\end{corollary}

\begin{proof}
By uncorrelatedness and
$\E\norm{A^{m}r_k}^2=\tr\bigl(A^{m}\Cov(r_k)(A^{m})^{*}\bigr)
 \le\opnorm{A^m}^2\tr\Cov(r_k)$,
\[
  \E\norm{(\Lacc r)_n}^2=\sum_{m=0}^{n-1}\E\norm{A^{m}r_{n-m}}^2
  \le\rho^2\sum_{m=0}^{n-1}\opnorm{A^m}^2 .
\]
Summing over $n$ gives $\E\norm{\Lacc r}^2\le\rho^2\tilde\Gamma_N^2$, and $\opnorm{\Sa_N}^2\le\beta$, $\opnorm{\Syn_N}\le\alpha_N^{-1/2}$ give the noiseless statement from \eqref{eq:exact-error}; the general case follows from $\norm{x+y}^2\le2\norm x^2+2\norm y^2$.
\end{proof}

\begin{corollary}\label{cor:worstcase}
Assume $\GG_N$ is a frame. If $\norm{r_k}\le\rho$ for all $k$, then with $a:=\opnorm A$,
$\Theta_n:=\sum_{m=0}^{n-1}a^{m}$ and
$\Gamma_N:=(\sum_{n\le N}\Theta_n^2)^{1/2}$,
\begin{equation}\label{eq:propagation}
  \norm{f-\hat f}\ \le\ \alpha_N^{-1/2}\bigl(\sqrt{\beta}\,\Gamma_N\rho+\norm{\varepsilon}\bigr).
\end{equation}
Taking $\rho=\bar\rho_N$ from Lemma \ref{lem:uniform} gives this bound with probability at least $\delta$, retaining the realised noise norm. In the noiseless case the second term vanishes. For a deterministic high-probability bound on both terms, choose failure probabilities $q_r,q_\varepsilon>0$ with $q_r+q_\varepsilon<1$, use confidence $1-q_r$ in Lemma \ref{lem:uniform}, and replace $\norm{\varepsilon}$ by
$\sqrt{\tr\Sigma_\varepsilon}
+\sqrt{2\opnorm{\Sigma_\varepsilon}\log(1/q_\varepsilon)}$.
The resulting confidence is at least $1-q_r-q_\varepsilon$, without any independence assumption between the two errors.
\end{corollary}

\begin{proof}
As above, with the triangle inequality
$\norm{(\Lacc r)_n}\le\sum_{m<n}\opnorm{A}^{m}\rho=\Theta_n\rho$ in place of
the quadrature step, and add
$\opnorm{\Syn_N}\norm{\varepsilon}$ using \eqref{eq:exact-error}.
\end{proof}

\begin{remark}\label{rem:relaxations}
The two corollaries use different information, quadrature requires temporal
uncorrelatedness of $r_k$, which need not hold for the posterior errors, while
the worst-case bound uses pathwise norms. They are not interchangeable under
an arbitrary correlated prior. If $\norm{A^m}=1$ for all $m$, the accumulation
factors scale as $\tilde\Gamma_N\sim N/\sqrt2$ and
$\Gamma_N\sim N^{3/2}/\sqrt3$. For $\norm{A^m}=a^m$ with $a<1$, their ratio
tends to $\sqrt{(1+a)/(1-a)}$. For a general operator,
replacing $\norm{A^m}$ by $\norm A^m$ can introduce further loss.

The exact identity also retains the orientation of the covariance relative to the sampled accumulation and dual-frame operators, as well as the noise cross-covariance. Comparisons must use the same error quantity (mean square, root mean square, or a stated high-probability norm bound) and verify the assumptions of each corollary. An RMS relaxation that does not require temporal uncorrelatedness follows instead by applying the triangle inequality in $L^2$ at each accumulation step. With $\tau_N:=\max_{1\le k\le N}(\tr\Cov(r_k))^{1/2}$,
\begin{equation}\label{eq:rms-relaxation}
  \bigl(\E\norm{f-\hat f}^2\bigr)^{1/2}
  \le \alpha_N^{-1/2}\left(
    \sqrt\beta\,\Gamma_N\tau_N+\sqrt{\tr\Sigma_\varepsilon}\right).
\end{equation}
Indeed, $(\E\norm{(\Lacc r)_n}^2)^{1/2}\le\Theta_n\tau_N$; summing these squared bounds and applying the triangle inequality to \eqref{eq:exact-error} gives \eqref{eq:rms-relaxation}. This is the RMS bound compared with the exact identity in Figure \ref{fig:relaxations}.
\end{remark}


Three quantities enter the coarse propagation bound, the nominal frame constant, the accumulation factor, and the forcing-error bound. The lower frame bound is non-decreasing for nested sampling systems, but need not saturate. The statistical uncertainty may decrease at a fixed time as more residuals are collected, while the uniform bound over the expanding window need not decrease. We first distinguish actual powers of $A$ from their norm-based majorants.

\begin{proposition}\label{prop:Gamma}
Put $a=\opnorm A$ and define the majorant
\[
  \widetilde\Gamma_N^{(a)}:=
    \Bigl(\sum_{n=1}^N\sum_{m=0}^{n-1}a^{2m}\Bigr)^{1/2}.
\]
Then $\sqrt N\le\tilde\Gamma_N\le\widetilde\Gamma_N^{(a)}$, and:
\begin{enumerate}
  \item if $a<1$, then $\sqrt N\le\Gamma_N\le\sqrt N/(1-a)$ and
        $\widetilde\Gamma_N^{(a)}\le\sqrt{N/(1-a^2)}$;
  \item if $a=1$, then $\Gamma_N\sim N^{3/2}/\sqrt3$ and
        $\widetilde\Gamma_N^{(a)}\sim N/\sqrt2$;
  \item if $a>1$, then $\Gamma_N$ and $\widetilde\Gamma_N^{(a)}$ are comparable to $a^N$, with constants depending on $a$.
\end{enumerate}
Both actual factors $\Gamma_N$ and $\tilde\Gamma_N$ diverge at least as $\sqrt N$. No matching linear or exponential lower growth for $\tilde\Gamma_N$ follows from $a=1$ or $a>1$ alone. For example, a nilpotent operator can have norm greater than one while only finitely many powers contribute. Equality with the majorant holds when $\norm{A^m}=a^m$ for all $m$.
\end{proposition}

\begin{proof}
The term $m=0$ gives the lower bound $\sqrt N$ for $\tilde\Gamma_N$, and $\opnorm{A^m}\le a^m$ gives its upper bound.
For $a\ne1$, $\Theta_n=(1-a^n)/(1-a)$, while $\Theta_n=n$ at $a=1$.
Summing $\Theta_n^2$ and $\sum_{m<n}a^{2m}$ proves the majorant estimates.
\end{proof}

\begin{corollary}\label{cor:Nstar}
Fix a depth $M$ and a consistent separable prior with fixed $C_w$. On admissible horizons with $\alpha_N>0$, define $E(N,\delta):=\alpha_N^{-1/2}\sqrt\beta\,\Gamma_N\bar\rho_N$, the forcing contribution to Corollary \ref{cor:worstcase}. 
Suppose $\alpha_N\uparrow\alpha_\infty<\infty$ and $\kappa_{\mathrm{eff},M}>0$. Then $E(N,\delta)\to\infty$, and it attains its minimum at a finite admissible horizon.
\end{corollary}

\begin{proof}
By Remark \ref{rem:floor}, $\bar\rho_N\ge\kappa_{\mathrm{eff},M}>0$, while $\alpha_N^{-1/2}\ge\alpha_\infty^{-1/2}$ and $\Gamma_N\ge\sqrt N$. Thus the bound diverges. A finite-valued sequence on the admissible integer horizons that diverges to infinity attains a minimum among finitely many horizons.
\end{proof}

\begin{remark}\label{rem:no-Nstar}
Corollary \ref{cor:Nstar} concerns the forcing contribution $E(N,\delta)$ to a coarse bound. It does not assert that the exact second moment has the same minimizer. In the contracted-rotation experiment of Section \ref{sec:exp-horizon}, the tested-grid minimum of the median bound occurs at $N=10$, whereas those of the median exact and observed RMS occur at $N=24$.

Nor does the exact error necessarily decrease or saturate. In the periodic-convolution experiment with isotropic forcing, its RMS decreases to about $0.1531$ at $N=10$ and then increases to about $0.2015$ at $N=24$, although $\opnorm A=1$ (Figure \ref{fig:relaxations}). The rank-three prior in the same sampling geometry instead yields a decreasing error over the tested frame horizons. These behaviours depend on the propagated covariance, the reconstruction operator and the shared-noise covariance; they cannot be inferred from $\alpha_N$ or the forcing floor alone.
\end{remark}

\section{Numerical experiments}\label{sec:numerics}

The experiments use real block rotations and complex periodic convolution. In the real case, $\Hi=\R^{12}$ and $p=2$. Each design draws six rotation angles independently and uniformly on $[0,2\pi)$, conjugates the resulting block rotation by an independent Haar orthogonal matrix, and takes the first two columns of another independent Haar orthogonal matrix as sensors.
Thus $\opnorm A=1$ and the base-bank Bessel bound is $\beta=1$. Signals are drawn uniformly on the unit sphere, so $\norm f=1$; for frame horizons the fixed-parameter reconstruction error is independent of $f$. Raw measurement noises are independent across times, sensor banks and channels, with variance $\sigma^2$. Reuse of those measurements in the residuals induces correlations, which are retained in the likelihood and in \eqref{eq:joint-covariance}. Complex forcing and noise are proper, standard complex innovations are $(X+iY)/\sqrt2$, with independent standard real Gaussian $X,Y$.

The amplitude convention is $\tr C_w=\varepsilon^2$, with $\varepsilon=0.1$ throughout. Unless stated otherwise, the spatial prior has rank three, with $C_w=\varepsilon^2WW^*/3$; the columns of $W$ are the first three left singular vectors spanning the depth-$M$ residual sensor bank, so $\ran C_w\subseteq S_{M-1}$. The temporal kernel is $K(k,\ell)=\rho^{|k-\ell|}$, with $\rho=0.9$. All Gaussian estimators are supplied with the true spatial covariance orientation, so the experiments do not assess estimation of an unknown spatial support.

Table \ref{tab:orbit} reports medians of per-design RMS errors over 25 designs, with 200 independent trajectories per design and design interquartile ranges (IQRs). The periodic table uses one deterministic sampling design and 400 trajectories per horizon. Its approximate $95\%$ Monte Carlo intervals are $[\sqrt{\max(\widehat m-1.96s,0)},\sqrt{\widehat m+1.96s}]$, where $\widehat m$ is the empirical MSE and $s$ its estimated standard error. They describe Monte Carlo uncertainty at fixed design, not design variability or simultaneous confidence bands. Figure sweeps use 200 trajectories per point.


\subsection{Transport with an unmodeled source}\label{sec:exp-transport}

This example takes an image transported by a known flow while an unmodeled
source injects a perturbation in a fixed region.  The goal is to recover the
initial image $f$, and with it the whole nominal trajectory $A^nf$, from
samples of $y_n=A^nf+\delta_n$ taken with few sensors.  A sensor at every
pixel would make the problem trivial: measuring at $n=0$ suffices, since
$\delta_0=0$.  Sensor count is limited by cost and by space, so few sensors
must be combined with several measurement times, and then the samples carry
the accumulated deviation $\delta_n=\sum_{k=1}^{n}A^{\,n-k}w_k$, which depends
on every earlier forcing.  The scheme of Section \ref{sec:double} addresses
this: the residuals between the sensors and their auxiliaries determine the
one-step forcing $w_n$, accumulating $\hat w$ through $A$ gives $\hat\delta_n$,
and $G^*\hat\delta_n$ is subtracted from the samples before the frame is
inverted.  The question is whether a sensor count below that of dense
sampling suffices to recover $f$.

\paragraph{Evolution operator.}
We use $L\times L$ images with $L=64$ and indices modulo $L$.  The operator
combines horizontal transport of velocity $s$ with isotropic diffusion of
intensity $\nu$.  Being a periodic convolution, it is defined through its
symbol in the discrete Fourier transform, $\widehat{Af}=\hat a\,\hat f$, with
\begin{equation}\label{eq:symbol-ad}
  \hat a(\xi_x,\xi_y)=e^{-is\xi_x}\,
    e^{-\nu\left[(2-2\cos\xi_x)+(2-2\cos\xi_y)\right]} .
\end{equation}
The first factor shifts the image $s$ pixels in the $x$ direction; the second
is the discrete heat kernel, which blurs it in both.  For integer $s$ and
$\nu=0$, $A$ is the circular translation $(Af)(x,y)=f(x-s,y)$.  For
non-integer $s$ the translation cannot be realised by moving pixels: the image
moves by a fraction of a pixel and the values on the grid must be interpolated.
The natural way to do this is on the Fourier side, multiplying by the phase
$e^{-is\xi_x}$, which is band-limited interpolation.\footnote{With $L$ even the
Nyquist frequency $\xi_x=\pi$ coincides with $-\pi$, so for $A$ to map real
images to real images its symbol must be real there.  Since $e^{-is\pi}$ is not
real for non-integer $s$, we use the nearest real value of modulus one,
$\operatorname{sign}(\cos\pi s)$, which is $-1$ for $s=0.75$.  The convention
affects that frequency alone.}

\paragraph{Sensors.}
The sensors are point evaluations, $g_i=e_{x_i}$, placed in every row and in
one of every $m=4$ columns, so that there are $L^2/m=1024$ of them.  At depth
$M$ each sensor is complemented by the auxiliaries $(A^*)^kg_i$,
$k=1,\dots,M$, for a total of $1024(M+1)$: at $M=2$ this is $3072$, against
the $4096$ of dense sampling.

By definition of the adjoint, the auxiliary measures today what the original
sensor will measure $k$ steps later, $\ip{y}{(A^*)^kg_i}=\ip{A^ky}{g_i}$.
Under integer translation without diffusion this is simply another pixel, the
one $ks$ positions upstream of $x_i$.  In the fractional case with diffusion,
$A^k$ moves the image by a fraction of a pixel and blurs it, so $(A^*)^kg_i$
is no longer a point evaluation: it is a weighted average of the pixels near
$x_i-ks$, whose width grows with $k$.  Assumption \ref{ass:realisable} then
requires that sensors with that response can be built.

\paragraph{Perturbation.}
The forcing follows the separable prior of Assumption \ref{ass:prior}, with
$K(n,n')=\rho^{|n-n'|}$ and $\rho=0.9$, that is, a vector $AR(1)$ process.
Spatially, $C_w=DKD$, where $K$ is a squared-exponential kernel of length
scale $\ell$ and $D$ is the diagonal of a Gaussian window of width five pixels
centred at the source, so that the forcing acts only on a fixed region of the
image.  We take $\ell=2$ (localized forcing) and $\ell=8$ (smooth forcing),
and the scale of $C_w$ is fixed so that $\norm{\delta_N}/\norm f$ equals
$0.05$ or $0.2$ in root mean square.  The forcing stays at the source, but the
accumulated deviation $\delta_n$ is transported and blurred by $A$, forming a
wake in the direction of the flow.  Measurement noise is independent Gaussian
with $\sigma=10^{-3}$, and the horizon is $N=8$.

\begin{figure}[ht]
  \centering
  \includegraphics[width=\linewidth]{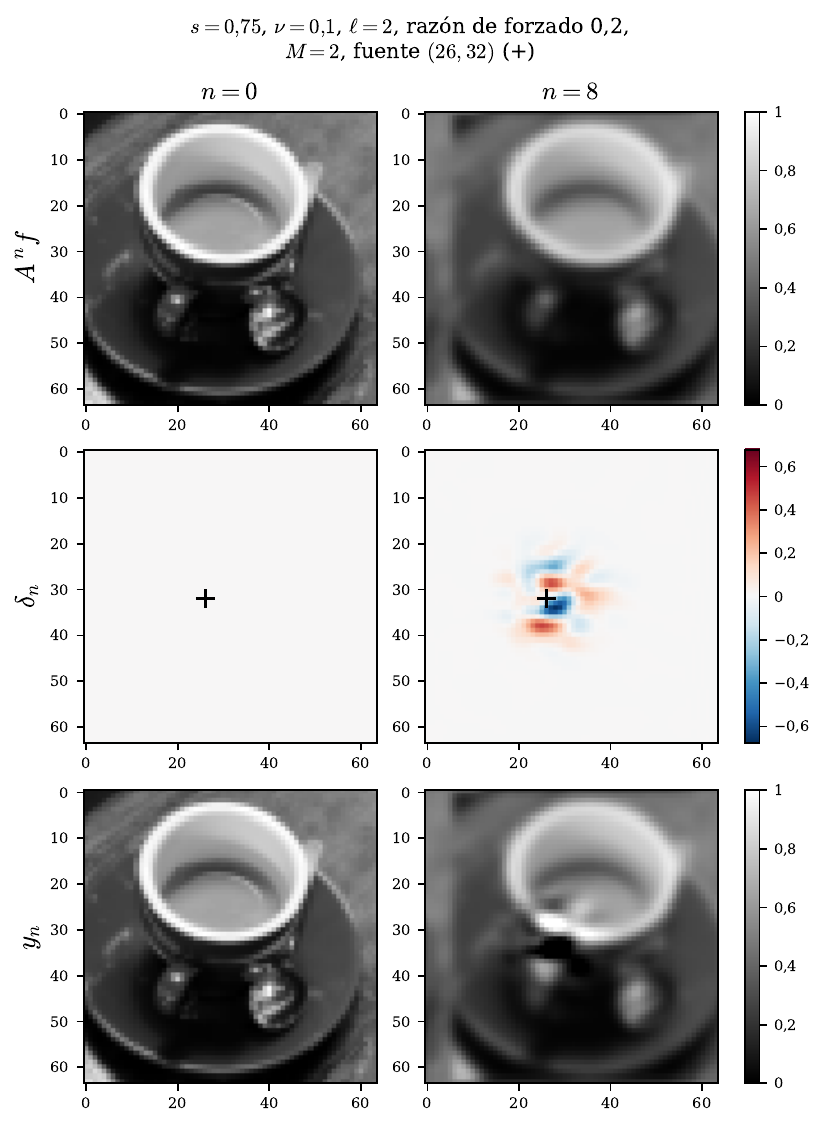}
  \caption{Trajectory of the advection--diffusion example ($s=0.75$,
  $\nu=0.1$, $\ell=2$, $\norm{\delta_N}/\norm f=0.2$) at $n=0$ and $n=8$.
  Top, the nominal trajectory $A^nf$, which is to be recovered; middle, the
  accumulated deviation $\delta_n$, generated by the source (cross) and
  dragged by the flow; bottom, $y_n=A^nf+\delta_n$, what the sensors actually
  measure.  In the top and bottom rows values outside $[0,1]$ are
  saturated.}\label{fig:ad-frames}
\end{figure}

\paragraph{Recoverability.}
Sampling takes inner products, $\ip{A^nf}{g_i}=\ip{f}{(A^*)^ng_i}$, so $f$ is
stably recoverable if and only if the orbit $\GG_N$ is a frame.  If it is not,
there are components of $f$ that no sensor sees at any time, and estimating
them would require prior information about the signal rather than about the
forcing.

Since the $y$ direction is not subsampled, a transform in $y$ separates the
problem into one-dimensional problems on $\mathbb C^{L}$, one per vertical
frequency, with sensors at $x_i=im$ and $L/m$ of them per row.  For a periodic
convolution sampled by point evaluations with step $m$, recoverability is
characterised in Remark \ref{rem:parity}: on each alias class
$\{q+(L/m)t\}_{t=0}^{m-1}$ the orbit analysis is a Vandermonde matrix with
nodes $\hat a(\xi_{q+(L/m)t})$, and $\GG_N$ is a frame if and only if
$N\ge m-1$ and those $m$ nodes are distinct in every class.  For an integer
translation this is equivalent to $\gcd(s,m)=1$: with $s=1$ and $m=4$, in four
instants every column of the image passes in front of a sensor.  If instead
$m\mid s$, every iterate stays on the original sublattice and no horizon gives
a frame.

\begin{remark}[Sensors with a footprint]\label{rem:footprint}
If the sensors are not point evaluations but $g_i=\varphi(\cdot-x_i)$ with
$\varphi$ a kernel averaging neighbouring pixels, the matrix of each alias
class is multiplied by the diagonal of the values
$\overline{\hat\varphi(\xi_{q+(L/m)t})}$.  Beyond the condition on the nodes
one then needs $\hat\varphi(\xi_k)\neq0$ at every grid frequency: a frequency
at which $\hat\varphi$ vanishes is invisible at any horizon.  For integer
translation the sample of sensor $i$ at time $n$ is $h(x_i-ns)$ with
$h(t)=\ip{f}{\varphi(\cdot-t)}$ the image seen through the kernel, and
$\hat h=\hat f\,\overline{\hat\varphi}$.  Numerically, a uniform average of
width $w$ fails to give a frame when $\gcd(w,L)>1$; so do B-splines at scale
$m$, which contain that factor, while a Gaussian of width one pixel gives a
frame with $\alpha_N\approx10^{-6}$.  With point sensors $\varphi=\delta$,
$\hat\varphi\equiv1$, the condition holds trivially, and that is the choice
made here.
\end{remark}

For the operator \eqref{eq:symbol-ad}, diffusion without transport ($s=0$,
$\nu>0$) has a real and even symbol and, by Remark \ref{rem:parity}, gives no
frame at any $N$; numerically $\alpha_N=0$ to machine precision.  Transport
breaks the parity and makes the operator admissible, but diffusion reduces the
modulus of the nodes and with it the frame constant: at $s=0.75$, $\alpha_8$
falls from $1.39$ at $\nu=0$ to $2.2\cdot10^{-2}$ at $\nu=0.1$.  Since
$\opnorm{\Syn_N}\le\alpha_N^{-1/2}$, a smaller frame constant means more
amplification of the measurement noise, and the oracle error, which comes from
noise alone, grows from $1.6\cdot10^{-3}$ to $5.8\cdot10^{-3}$.  The frame
constant moreover saturates shortly after $N=m-1$, so extending the horizon
does not recover that loss.  The choice of $s$ relative to $m$ also matters:
at $s=1.3$ two nodes of one class nearly coincide, because
$\tfrac32\pi s\approx2\pi$, and $\alpha_8$ is almost ten times smaller than at
$s=0.75$.

\paragraph{The cost of depth.}
The auxiliary sensor of order $k$ sits $ks$ positions upstream.  With integer
translation and point sensors, $(A^*)^kg_i=e_{x_i-ks}$ coincides with an
already installed sensor exactly when $m\mid ks$, and if $\gcd(s,m)=1$ this
means $m\mid k$: for $k=1,\dots,m-1$ the auxiliaries occupy new columns.  The
favourable case $m\mid s$, in which the installed array supplies the auxiliary
samples at no extra cost, is precisely the one that gives no frame.  At
$s=0.75$ the displacement is not an integer and the auxiliaries do not
coincide with existing sensors either.  Depth is therefore paid for here, each
level adding $1024$ sensors, and the relevant question is not whether the
auxiliaries are free but whether a total below that of dense sampling
suffices.

\paragraph{Results.}
We compare four reconstructions: the one that ignores the forcing and uses
only the $1024$ base sensors (\emph{uncorrected}); the one that uses $3072$
point sensors, the same count as the scheme at $M=2$, and inverts the frame
directly without estimating the forcing (\emph{naive}); the Gaussian scheme of
Definition \ref{def:estimator}; and the \emph{oracle}, which subtracts the
true deviation and whose error comes from noise alone.  Dense sampling, with
$4096$ sensors at $n=0$, is included for reference.  Entries are medians of
the relative error $\norm{\hat f-f}/\norm f$ over $200$ realisations.

\begin{table}[ht]
  \centering\small
  \caption{Reconstruction error in the advection--diffusion example
  ($s=0.75$, $\nu=0.1$, $M=2$, $N=8$).  Medians over $200$ realisations;
  sensor counts in brackets.}\label{tab:transport}
  \begin{tabular}{cccccccc}
    \toprule
    $\ell$ & $\norm{\delta_N}/\norm f$ & uncorrected & naive & Gaussian
      & oracle & dense \\
     & & ($1024$) & ($3072$) & ($3072$) & & ($4096$) \\
    \midrule
    2 & 0.05 & 0.0395 & 0.0293 & 0.0065 & 0.0058 & 0.0023 \\
    2 & 0.2  & 0.1563 & 0.1169 & 0.0096 & 0.0058 & 0.0023 \\
    8 & 0.05 & 0.0273 & 0.0218 & 0.0058 & 0.0058 & 0.0023 \\
    8 & 0.2  & 0.1068 & 0.0867 & 0.0059 & 0.0058 & 0.0023 \\
    \bottomrule
  \end{tabular}
\end{table}

Table \ref{tab:transport} shows three effects.  First, at equal sensor count
the Gaussian scheme reduces the error of the naive reconstruction by a factor
between $3.8$ and $15$, whereas tripling the sensors without a model of the
forcing barely reduces it, from $0.1563$ to $0.1169$ at $\ell=2$ and ratio
$0.2$.  The error of the naive reconstruction is dominated by the uncorrected
forcing, not by a shortage of data.  Second, with smooth forcing ($\ell=8$)
the Gaussian scheme attains the oracle, although the residuals at $M=2$
observe only half the columns: by Lemma \ref{lem:lift}, correlation in the
prior predicts the invisible part of the forcing from the visible one.  This
is the situation of Remark \ref{rem:rank}, in which the rank condition fails
and the reconstruction is nevertheless as good as if the forcing were known.
With localized forcing ($\ell=2$) that correlation is weak and an excess
remains, small at ratio $0.05$ and larger at ratio $0.2$; this is
$\kappa_{\mathrm{eff},M}>0$ appearing in the reconstruction error.  Third,
dense sampling lies below the oracle.  The oracle error is noise amplified by
$\alpha_N^{-1/2}$ and diffusion reduces $\alpha_N$, while dense sampling,
which measures at $n=0$, does not suffer it.  The correction of the forcing
does not depend on the operator, but the error floor does: with pure transport
($s=1$, $\nu=0$) the oracle falls below dense sampling,
$1.54\cdot10^{-3}$ against $2.28\cdot10^{-3}$, because each pixel is measured
more than once along the horizon.

\begin{figure}[ht]
  \centering
  \includegraphics[width=0.8\linewidth]{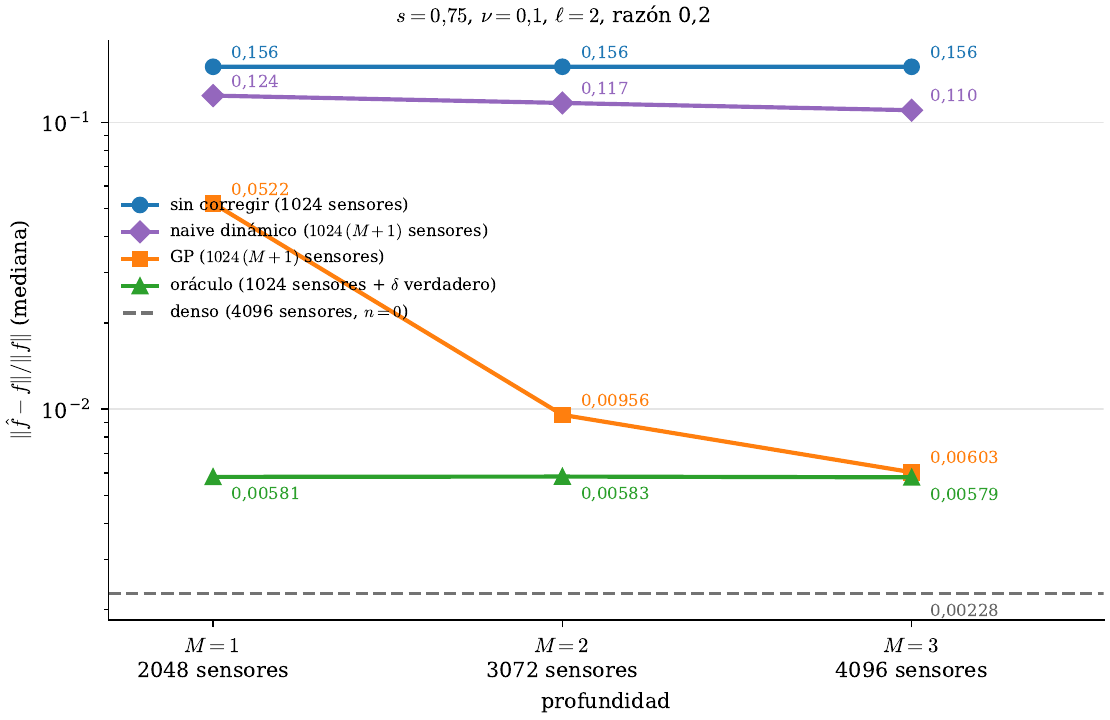}
  \caption{Relative reconstruction error against augmentation depth in the
  advection--diffusion example ($s=0.75$, $\nu=0.1$, $\ell=2$,
  $\norm{\delta_N}/\norm f=0.2$), with the sensor count of each scheme.
  Logarithmic scale.  The Gaussian error falls from $0.0522$ at $M=1$ to
  $0.0060$ at $M=3$, approaching the oracle, while the same sensors spent
  without a model of the forcing move the naive error only from $0.124$ to
  $0.110$.}\label{fig:ad-depth}
\end{figure}

Figure \ref{fig:ad-depth} separates the two uses of a sensor.  Along the
Gaussian curve the depth enlarges $S_{M-1}$, so $\kappa_{\mathrm{eff},M}$
decreases and the error approaches the oracle; along the naive curve the same
measurements enter the frame inversion instead, and the error hardly moves.
The gap between the two curves at a fixed abscissa is the value of modelling
the forcing at a fixed budget.

\begin{figure}[ht]
  \centering
  \includegraphics[width=\linewidth]{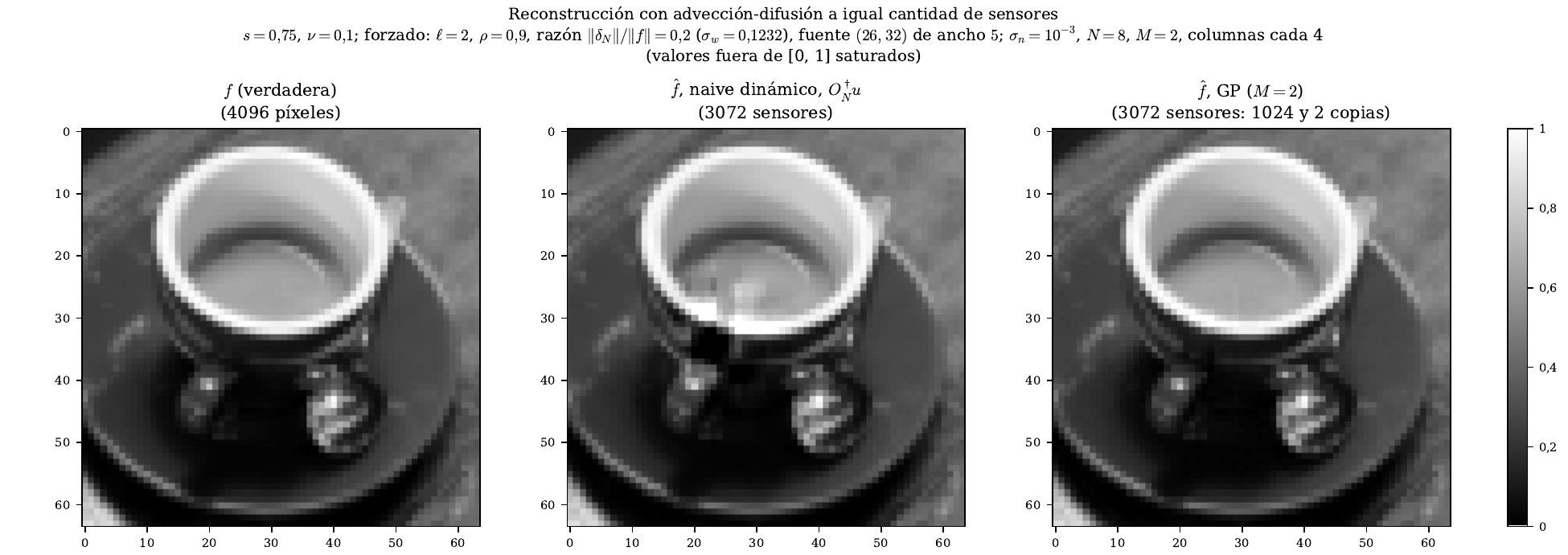}
  \caption{Reconstructions of one realisation, same configuration as
  Figure \ref{fig:ad-depth}: the initial image $f$, the naive reconstruction
  from $3072$ sensors, and the Gaussian reconstruction from the same $3072$.
  The artefact left by the uncorrected source is visible in the
  former.}\label{fig:ad-reconstruction}
\end{figure}

\paragraph{Design at a fixed budget.}
A given number of sensors can be divided between base sensors and auxiliaries
in several ways: fewer columns with more depth, or more columns with less
depth.  Table \ref{tab:budget} compares those configurations.

\begin{table}[ht]
  \centering\small
  \setlength{\tabcolsep}{4pt}
  \caption{Reconstruction error at a fixed budget $P$ in the
  advection--diffusion example ($s=0.75$, $\nu=0.1$, ratio $0.2$, $N=8$).
  Medians over $200$ realisations.  The configuration $P=3072$, $m=2$ is
  omitted because $P/p$ is not an integer.}\label{tab:budget}
  \begin{tabular}{cccccccccc}
    \toprule
    $P$ & $m$ & $p$ & $M$ & $\alpha_8$ & $\ell$ & uncorrected & naive
      & Gaussian & oracle \\
    \midrule
    2048 & 2 & 2048 & 0 & $3.09\cdot10^{-1}$ & $2$ & $0.1240$ & $0.1240$ & --- & --- \\
    2048 & 2 & 2048 & 0 & $3.09\cdot10^{-1}$ & $8$ & $0.0912$ & $0.0912$ & --- & --- \\
    2048 & 4 & 1024 & 1 & $2.17\cdot10^{-2}$ & $2$ & $0.1563$ & $0.1240$ & $0.0522$ & $0.0058$ \\
    2048 & 4 & 1024 & 1 & $2.17\cdot10^{-2}$ & $8$ & $0.1068$ & $0.0912$ & $0.0059$ & $0.0058$ \\
    2048 & 8 &  512 & 3 & $4.91\cdot10^{-5}$ & $2$ & $0.2193$ & $0.1240$ & $0.0805$ & $0.0616$ \\
    2048 & 8 &  512 & 3 & $4.91\cdot10^{-5}$ & $8$ & $0.1496$ & $0.0912$ & $0.0616$ & $0.0616$ \\
    \midrule
    3072 & 4 & 1024 & 2 & $2.17\cdot10^{-2}$ & $2$ & $0.1563$ & $0.1169$ & $0.0096$ & $0.0058$ \\
    3072 & 4 & 1024 & 2 & $2.17\cdot10^{-2}$ & $8$ & $0.1068$ & $0.0867$ & $0.0059$ & $0.0058$ \\
    3072 & 8 &  512 & 5 & $4.91\cdot10^{-5}$ & $2$ & $0.2203$ & $0.1169$ & $0.0630$ & $0.0620$ \\
    3072 & 8 &  512 & 5 & $4.91\cdot10^{-5}$ & $8$ & $0.1500$ & $0.0867$ & $0.0620$ & $0.0620$ \\
    \midrule
    \multicolumn{6}{l}{dense ($4096$ sensors, $n=0$)}
      & \multicolumn{4}{c}{$0.0023$} \\
    \bottomrule
  \end{tabular}
\end{table}

The budget is best spent on base sensors, and the reason is the separation of
Section \ref{sec:lift-sub}: depth buys visibility of the forcing, but only the
base bank buys nominal stability, and the two are not interchangeable.  Going
from $m=4$ to $m=8$ at a fixed budget trades one for the other and loses:
$\alpha_8$ collapses from $2.17\cdot10^{-2}$ to $4.91\cdot10^{-5}$, and with
it the oracle floor rises from $0.0058$ to $0.0616$, so that no amount of
depth can help.  At $P=2048$ the extreme in the other direction, $m=2$ with
$M=0$, has the best frame constant of the three, $3.09\cdot10^{-1}$, but no
residuals at all and therefore no correction: it is stuck at $0.1240$.  The
useful configurations are interior, which is the quantitative form of the
statement that the two geometric requirements must be met separately.

In summary, the example shows three behaviours.  At equal sensor count,
correction through residuals and the Gaussian posterior removes nearly all of
the contribution of the forcing, while adding sensors without a model of it
barely improves the reconstruction.  The depth required depends on the
structure of the forcing: a localized source calls for more observed
directions, whereas for a smooth one the prior compensates what the residuals
do not see.  Finally, the dynamics fixes the error floor: transport preserves
the information along the horizon, whereas diffusion attenuates it and limits
what dynamical sampling can gain over dense sampling.

\subsection{Depth is not a substitute for the orbit}\label{sec:exp-orbit}

We fix $M=2$. A single time-zero snapshot then uses only the six vectors of $\GG_2$ in a twelve-dimensional state space, and cannot determine an arbitrary $f$. The rank-three forcing is supported in $S_1=\operatorname{span}\{g_j,A^*g_j\}$, so $\kappa_M=0$ and the
visibility condition \eqref{eq:visibility} holds. We use $\rho=0.9$ and
$\sigma=10^{-3}$. For the sampled designs the base-bank orbit first becomes
a frame at $N=5$. Its median lower frame bound is $6.60\cdot10^{-4}$ at
that horizon, increasing to $0.0444$ at $N=8$ and $0.305$ at $N=15$.

\begin{table}[ht]
\centering\small
\caption{Initial-state reconstruction at depth $M=2$. Entries are medians of per-design RMS errors, with the design interquartile range in brackets; 25 designs and 200 trajectories per design. The $N=3$ row uses minimum-norm least squares and is outside the frame theorem.}\label{tab:orbit}
\begin{tabular}{rccc}
\toprule
$N$ & median $\alpha_N$ & uncorrected RMS & corrected RMS\\
\midrule
3 & 0 & 0.6375 [0.6162, 0.6682] & 0.5791 [0.5733, 0.5861]\\
5 & 0.00066 & 2.0415 [0.7804, 3.3251] & 0.0373 [0.0177, 0.1205]\\
8 & 0.0444 & 0.3980 [0.2279, 0.5914] & 0.0082 [0.0053, 0.0172]\\
11 & 0.152 & 0.3142 [0.1763, 0.4772] & 0.0062 [0.0049, 0.0076]\\
15 & 0.305 & 0.2905 [0.1791, 0.4972] & 0.0057 [0.0046, 0.0065]\\
\bottomrule
\end{tabular}
\end{table}

For $N=3$, both reconstructions use minimum-norm least squares; this rank-deficient row is outside Theorem \ref{thm:main}. It has a signal-dependent projection error as well as forcing and measurement errors. At $N=5$ the median uncorrected RMS is $2.0415$, worse than the zero estimator whose error is one. Correction reduces it to $0.0373$, a ratio of medians of approximately $55$. At $N=8$ the corresponding ratio is approximately $49$. The IQRs show substantial variation between designs near the frame threshold. Correction remains beneficial at the longer tested horizons, but these numerical ratios are specific to the sampling and prior designs.

\begin{remark}\label{rem:two-scales}
Visibility of the forcing and recoverability of the initial state concern different objects. The condition $\ran C_w\subseteq S_{M-1}$ makes the forcing directly visible from depth-$M$ residuals; it does not make a time-zero snapshot a frame for the full state space. The orbit horizon must still supply nominal observability. In this experiment, depth $M=2$ is sufficient for forcing visibility, while the base-bank orbit first spans the state space at $N=5$ and becomes better conditioned thereafter.
\end{remark}

\subsection{Convolution on a sublattice}\label{sec:exp-lattice}

A finite sensor set cannot form a finite-horizon dynamical frame for an infinite-dimensional state space. The canonical infinite-lattice setting therefore uses countably many sensors on $m\Z$. Here we use the finite periodic state space $\Hi=\ell^2(\Z_{64})$, with $m=8$, $p=8$ sensors $g_j=\chi_{\{jm\}}$, and convolution $A=\mathcal F^*\operatorname{diag}(\hat a_k)\mathcal F$, where $\mathcal F$ is the unitary discrete Fourier transform and
\[
  \hat a_k=e^{1.3i\xi_k},\qquad \xi_k=\frac{2\pi\tilde k}{64},\qquad
  \tilde k=\begin{cases}k,&0\le k<32,\\k-64,&32\le k<64.\end{cases}
\]
This balanced branch specifies a finite discrete symbol; the noninteger phase is not being used as a continuous $2\pi$-periodic symbol. The operator is unitary, so $\opnorm A=1$. The forcing has rank three inside $S_1$, with $M=2$, $\rho=0.9$, $\varepsilon=0.1$ and $\sigma=10^{-3}$.

The base sensors span eight of the 64 dimensions. The depth-two time-zero bank has at most 24 independent vectors and cannot determine arbitrary $f$. The orbit first becomes a frame at $N=7=m-1$, where $\alpha_7=0.00783013$ for the specified Fourier branch.

\begin{table}[ht]
\centering\small
\caption{Periodic convolution on the balanced Fourier grid, $L=64$, $m=8$, $\theta=1.3$, $M=2$, rank-three forcing, $\rho=0.9$, amplitude $0.1$, and $\sigma=10^{-3}$. Exact and observed columns are RMS; brackets give approximate 95\% Monte Carlo intervals from 400 proper complex trajectories.}\label{tab:lattice}
\begin{tabular}{rcccc}
\toprule
$N$ & $\alpha_N$ & exact RMS & observed RMS [95\% MC interval] & uncorrected RMS\\
\midrule
7 & 0.007830 & 0.03324 & 0.03327 [0.03273, 0.03380] & 0.35808\\
8 & 0.026040 & 0.02007 & 0.02006 [0.01979, 0.02033] & 0.20968\\
12 & 0.175212 & 0.01073 & 0.01065 [0.01055, 0.01076] & 0.15138\\
20 & 0.748908 & 0.00759 & 0.00755 [0.00749, 0.00762] & 0.14772\\
24 & 1.525714 & 0.00703 & 0.00703 [0.00697, 0.00708] & 0.14944\\
\bottomrule
\end{tabular}
\end{table}

\begin{figure}[h]
\centering
\includegraphics[width=0.78\linewidth]{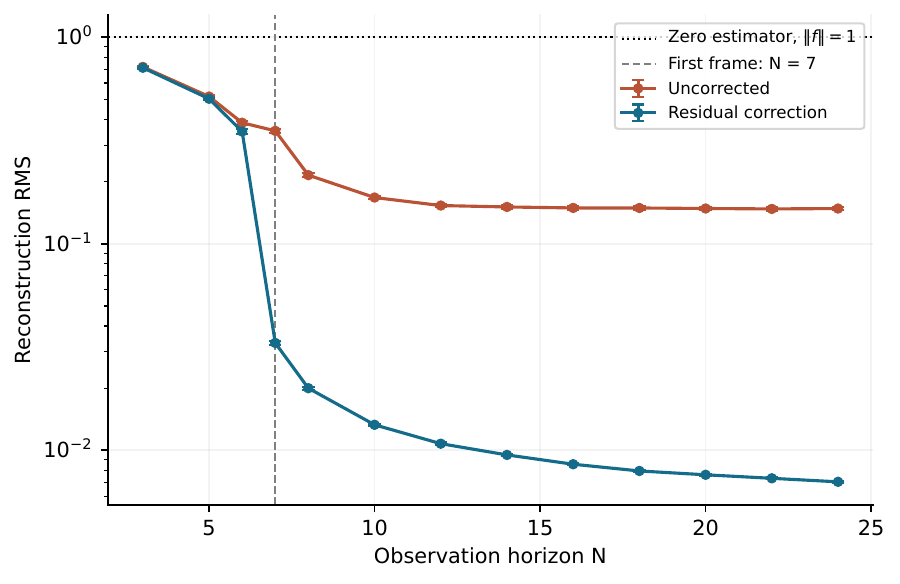}
\caption{Periodic-convolution reconstruction with and without residual correction at depth $M=2$. Points use 200 trajectories and error bars are approximate $95\%$ Monte Carlo intervals for RMS. The first frame horizon is $N=7$; earlier points use minimum-norm least squares and are outside Theorem \ref{thm:main}. The zero estimator has error one.}\label{fig:horizon}
\end{figure}

The exact RMS is $(\tr\Sigma_f)^{1/2}$, computed with the full joint covariance in \eqref{eq:joint-covariance}. All five predictions lie within the approximate Monte Carlo intervals. Independent calculations from the raw forcing and noise maps agree with the joint-block covariance to relative error of order $10^{-15}$. Omitting the cross-covariance raises the predicted RMS by approximately $10\%$ at $N=7$ and $18\%$ at $N=24$. The observed correction gain ranges from $10.5$ to $21.3$ across the five rows; it is $10.8$ at $N=7$ and $21.3$ at $N=24$. Figure \ref{fig:horizon} also includes horizons below the frame threshold, using minimum-norm least squares there.

\begin{remark}[Aliasing in the periodic model]\label{rem:parity}
Let $\Hi=\ell^2(\Z_L)$, $m\mid L$, $p=L/m$, with sensors at $jm$. On the Fourier side, the alias class indexed by $q\in\{0,\ldots,p-1\}$ is $\{q+pt:t=0,\ldots,m-1\}$. The orbit analysis on this class is a Vandermonde matrix with nodes $\hat a_{q+pt}$ and rows $n=0,\ldots,N$. Consequently, $\GG_N$ is a frame precisely when $N\ge m-1$ and the $m$ nodes are pairwise distinct in every alias class.

Even discrete symbols, $\hat a_k=\hat a_{-k\bmod L}$, can obstruct this condition. Reflection maps the class of $q$ to that of $-q\bmod p$; a class is invariant exactly when $2q\equiv0\pmod p$. In the present $L=64,m=p=8$ model, class $q=0$ contains the distinct indices $8$ and $56$, which have equal symbols under evenness. Thus even symbols fail the frame condition for this periodic geometry. This is a statement about this finite sampling model, not a universal assertion that every even convolution symbol fails for every sublattice. In an infinite lattice model, stable recovery requires a uniform lower bound on the corresponding fiber matrices; distinct nodes almost everywhere alone ensure neither that bound nor a finite-horizon frame.
\end{remark}


\subsection{Exact covariance, coarse bounds and conditional uncertainty}
\label{sec:exp-bounds}

\begin{figure}[h]
\centering
\includegraphics[width=\linewidth]{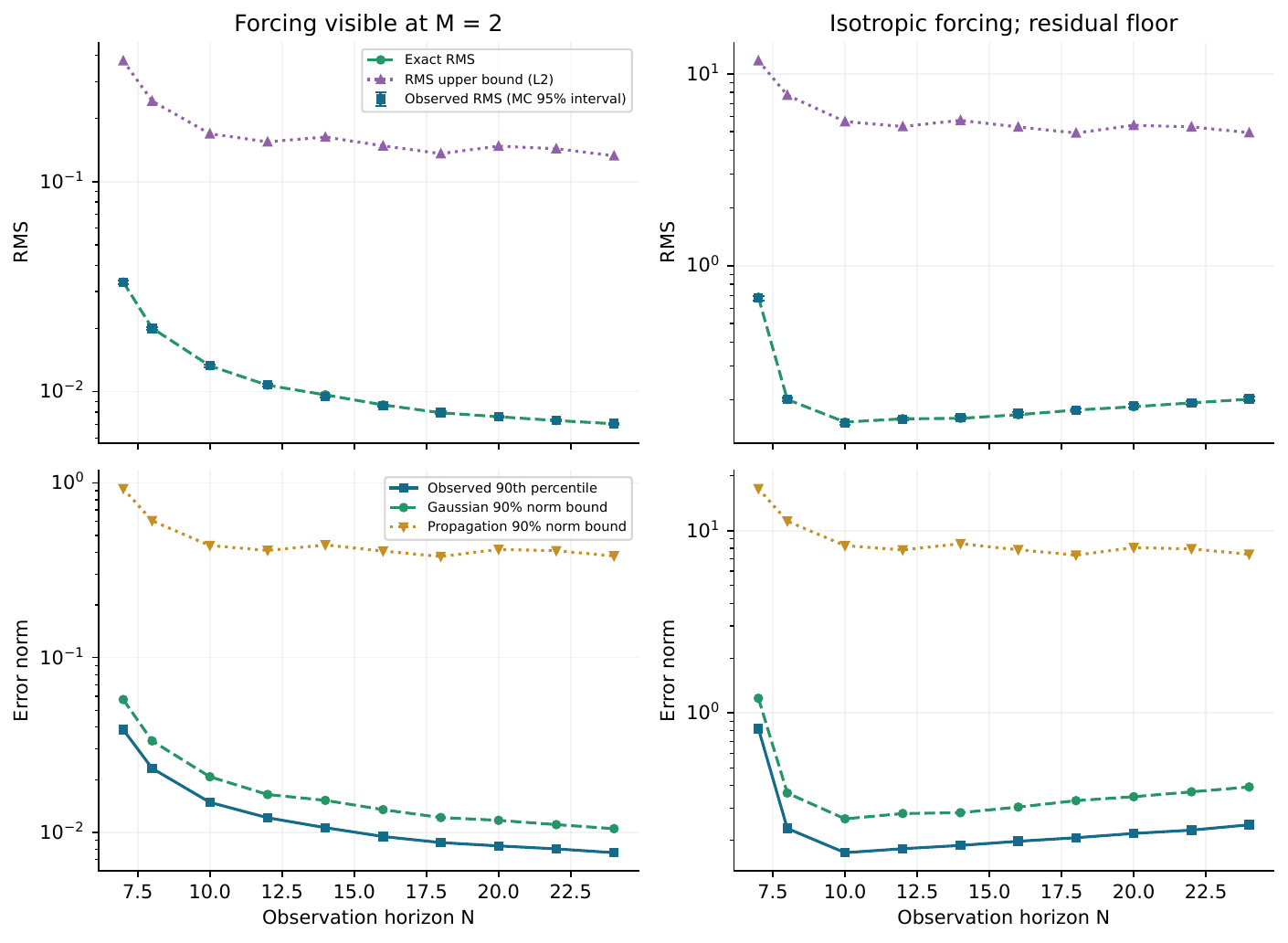}
\caption{Comparable metrics for the rank-three and isotropic priors in the periodic design, with $M=2$, $\rho=0.9$, $\sigma=10^{-3}$ and 200 trajectories per point. Top: exact and observed RMS against the valid $L^2$ bound \eqref{eq:rms-relaxation}. Bottom: empirical 90th percentiles against norm bounds with confidence at least $90\%$.
The propagation bound includes both forcing and measurement noise; the quadrature corollary is not used.}\label{fig:relaxations}
\end{figure}

Figure \ref{fig:relaxations} compares the rank-three prior with the isotropic prior $C_w=(0.1^2/64)I$ in the same periodic sampling geometry. The upper panels compare RMS with \eqref{eq:rms-relaxation}, the lower panels compare empirical 90th-percentile error norms with the Gaussian $90\%$ bound \eqref{eq:exact-hp} and the propagation bound of Corollary \ref{cor:worstcase}. The latter includes measurement noise and uses failure probabilities $q_r=q_\varepsilon=0.05$, giving total confidence at least $90\%$. All coarse bounds use the base-bank constant $\beta=1$.

The quadrature corollary is not applied, the posterior errors are temporally correlated. In Table \ref{tab:lattice}, the fraction $\norm{\Sigma_r-\operatorname{blockdiag}(\Sigma_r)}_F/\norm{\Sigma_r}_F$ is approximately $0.33$--$0.35$, so its additional hypothesis fails. The large gap between valid propagation bounds and exact covariance predictions illustrates the loss of orientation and covariance information in the coarse estimates.

For the isotropic prior the exact RMS attains its smallest tested value at $N=10$ (approximately $0.1531$), and then rises to $0.2015$ at $N=24$. Thus error growth does not require $\opnorm A>1$. These finite-grid observations do not imply a universal optimal horizon, and the forcing floor alone does not determine the reconstruction error.

\begin{figure}[h]
\centering
\includegraphics[width=0.78\linewidth]{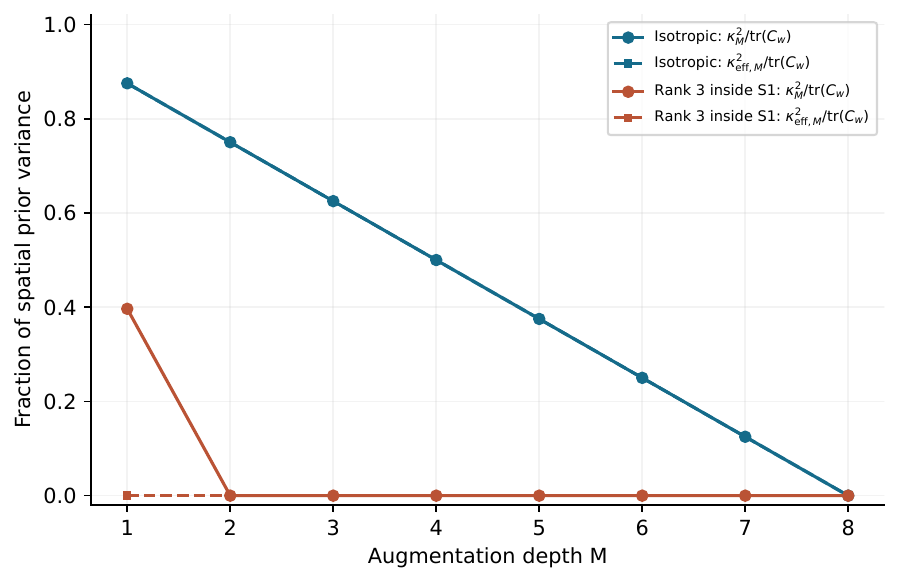}
\caption{Invisible prior variance $\kappa_M^2$ and conditional irreducible variance $\kappa_{\mathrm{eff},M}^2$, normalised by $\tr C_w$, against augmentation depth in the periodic design. Each spatial prior is held fixed as depth varies. For the rank-three prior, the conditional floor already vanishes at $M=1$, although invisible prior mass remains. These are forcing-uncertainty quantities, not initial-state reconstruction errors.}\label{fig:depth}
\end{figure}

Figure \ref{fig:depth} separates invisible prior mass from conditional uncertainty. The rank-three prior is held fixed inside $S_1$ while $M$ varies. At $M=1$ its invisible-mass fraction is approximately $0.3965$, yet the conditional floor is zero to numerical precision, the known spatial support allows the visible coordinates to determine the forcing. For the isotropic prior, the two fractions agree and decrease with depth.
Failure of the sufficient visibility condition alone therefore does not imply positive irreducible conditional uncertainty, in agreement with Example \ref{ex:predictability} and Remark \ref{rem:floor}.

\subsection{A bound minimizer and an error minimizer}\label{sec:exp-horizon}

To illustrate Corollary \ref{cor:Nstar}, we use the real rotation design contracted by $0.85$, so $\opnorm A=0.85$, with $M=1$, isotropic $C_w=(0.1^2/12)I$, $\rho=0.9$ and $\sigma=10^{-3}$. The frame bounds converge because the iterated sensors decay.
The conditional forcing floor is $\kappa_{\mathrm{eff},1}=0.1\sqrt{10/12}\approx0.0913$. We use 25 designs with 200 trajectories per design and the grid $N\in\{6,8,10,14,18,24,32\}$.

\begin{figure}[h]
\centering
\includegraphics[width=\linewidth]{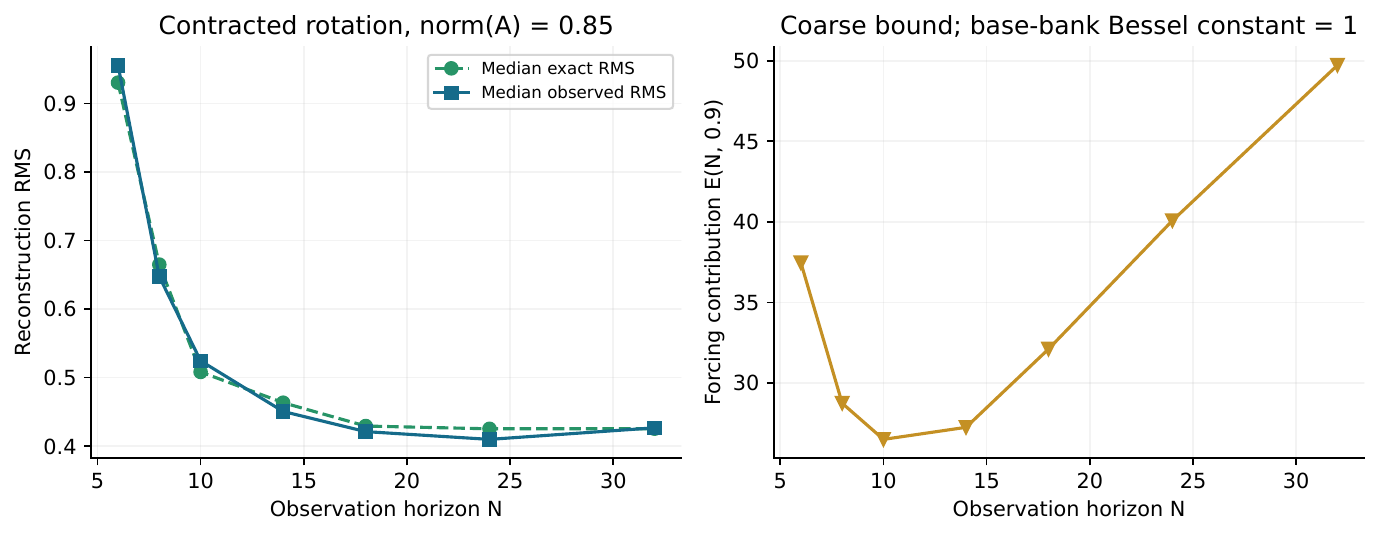}
\caption{Contracted-rotation experiment with $\opnorm A=0.85$, isotropic forcing and $M=1$. Left: medians of design-specific exact and observed RMS, with 25 designs and 200 trajectories per design.
Right: the median forcing contribution $E(N,0.9)$ computed using $\beta=1$. Its tested-grid minimizer differs from that of the exact RMS. The two panels show different quantities and use separate scales.}
\label{fig:horizon-tradeoff}
\end{figure}

The median forcing contribution $E(N,0.9)$ has its grid minimum $26.4926$ at $N=10$. In contrast, the median exact RMS has its grid minimum $0.42523$ at $N=24$, where the median observed RMS is $0.40981$; at $N=32$ the corresponding values are $0.42542$ and $0.42669$. Figure \ref{fig:horizon-tradeoff} displays this distinction.
The exact values at $N=24$ and $N=32$ are close, and no statistically resolved empirical preference between them is claimed. The grid minima of these design aggregates are neither design-wise nor global optima.
Moreover, $E$ is the forcing contribution only, rather than a bound on the entire noisy reconstruction error.

\section{Conclusions and outlook}\label{sec:outlook}

We have studied initial-state recovery from dynamical samples in the presence of unknown additive forcing. Augmenting the sensor system with its images under the adjoint nominal evolution produces residuals that cancel the initial-state contribution exactly. This cancellation permits forcing estimation to be carried out separately from initial-state reconstruction.
The sampling geometry remains essential at both stages, the augmented measurements determine which forcing components are observed, while the dynamical frame controls the stability of reconstruction from the corrected orbit samples.

For a centered Gaussian forcing model with fixed covariance parameters, this construction yields an unbiased reconstruction and an exact unconditional second-moment identity. The identity retains the covariance structure of the forcing-estimation error and its correlation with the measurement noise, including correlations induced by shared measurements. It therefore provides a basis for assessing sampling designs without discarding the orientation of the error relative to the reconstruction operator. The conditional covariance also identifies the forcing uncertainty that remains after observing all available residuals. This uncertainty need not vanish when the horizon increases, although its effect on initial-state reconstruction depends on how it propagates into the sampled data. The finite-dimensional and periodic-convolution experiments demonstrate the correction procedure and the different roles of sampling depth and observation horizon. For the specified periodic design, the exact RMS agrees with Monte Carlo at every reported frame horizon. The observed reduction from correction ranges from $10.5$ to $21.3$ across the table. The isotropic-prior experiment also shows that exact reconstruction error can increase with the horizon even for a unitary nominal operator, while the contracted example separates the minimizer of a coarse bound from that of the exact error.

Several questions remain beyond the present analysis. Estimating covariance parameters from the same residuals generally produces a nonlinear estimator, so extending the reconstruction guarantees to empirical Bayes requires accounting for parameter uncertainty. The experiments reported here supply the spatial covariance orientation and do not learn an unknown forcing subspace.
Larger trajectory samples and tests with misspecified spatial covariance are needed to assess that additional uncertainty. A further extension concerns state-dependent forcing, as generated by an uncertain evolution operator, for which the exogenous Gaussian model used here need not apply.
Non-separable covariance models would also allow temporal observations to convey information about spatial directions that the separable model leaves unresolved. Finally, the exact reconstruction covariance suggests a design problem in which sensor placement, augmentation depth and observation horizon are selected jointly under a measurement budget. Developing these extensions would broaden the applicability of residual-based correction while preserving the distinction between nominal observability and uncertainty in the forcing.

\section*{Declaration of generative AI technologies in the manuscript preparation process}
The authors used AI tools to assist with language revision and  manuscript organization. AI assistance was also used in reviewing the Python code for the numerical experiments. The authors reviewed and edited the assisted material and retain full responsibility for the numerical results, and final manuscript.

\bibliography{dynsampling_gp}

@article{DarouachZasadzinskiXu1994,
  author  = {Darouach, Mohamed and Zasadzinski, Michel and Xu, Shi Jie},
  title   = {Full-order observers for linear systems with unknown inputs},
  journal = {IEEE Trans. Automat. Control},
  volume  = {39},
  number  = {3},
  pages   = {606--609},
  year    = {1994},
}

@article{DarouachZasadzinski1997,
  author  = {Darouach, Mohamed and Zasadzinski, Michel},
  title   = {Unbiased minimum variance estimation for systems with unknown exogenous inputs},
  journal = {Automatica},
  volume  = {33},
  number  = {4},
  pages   = {717--719},
  year    = {1997},
}

@article{HouMuller1992,
  author  = {Hou, Min and M\"uller, Peter C.},
  title   = {Design of observers for linear systems with unknown inputs},
  journal = {IEEE Trans. Automat. Control},
  volume  = {37},
  number  = {6},
  pages   = {871--875},
  year    = {1992},
}

@article{HouPatton1998,
  author  = {Hou, Min and Patton, Ron J.},
  title   = {Input observability and input reconstruction},
  journal = {Automatica},
  volume  = {34},
  number  = {6},
  pages   = {789--794},
  year    = {1998},
}

@article{Hsieh2000,
  author  = {Hsieh, Chien-Shu},
  title   = {Robust two-stage {K}alman filters for systems with unknown inputs},
  journal = {IEEE Trans. Automat. Control},
  volume  = {45},
  number  = {12},
  pages   = {2374--2378},
  year    = {2000},
}

@article{GillijnsDeMoor2007,
  author  = {Gillijns, Steven and De Moor, Bart},
  title   = {Unbiased minimum-variance input and state estimation for linear discrete-time systems},
  journal = {Automatica},
  volume  = {43},
  number  = {1},
  pages   = {111--116},
  year    = {2007},
}

@article{Hautus1983,
  author  = {Hautus, Malo L. J.},
  title   = {Strong detectability and observers},
  journal = {Linear Algebra Appl.},
  volume  = {50},
  pages   = {353--368},
  year    = {1983},
}

@article{SundaramHadjicostis2007,
  author  = {Sundaram, Shreyas and Hadjicostis, Christoforos N.},
  title   = {Delayed observers for linear systems with unknown inputs},
  journal = {IEEE Trans. Automat. Control},
  volume  = {52},
  number  = {2},
  pages   = {334--339},
  year    = {2007},
}

@article{Kitanidis1987,
  author  = {Kitanidis, Peter K.},
  title   = {Unbiased minimum-variance linear state estimation},
  journal = {Automatica},
  volume  = {23},
  number  = {6},
  pages   = {775--778},
  year    = {1987},
}

@article{ACKM2026survey,
  author  = {Aldroubi, Akram and Cabrelli, Carlos and Krishtal, Ilya and Molter, Ursula},
  title   = {Dynamical Sampling: A Survey},
  journal = {La Matematica},
  volume  = {5},
  pages   = {},
  year    = {2026},
  doi     = {10.1007/s44007-026-00215-y},
  eprint  = {2511.10769},
  archivePrefix = {arXiv},
  primaryClass  = {math.FA}
}

@article{ACMT2017,
  author  = {Aldroubi, Akram and Cabrelli, Carlos and Molter, Ursula and Tang, Sui},
  title   = {Dynamical sampling},
  journal = {Applied and Computational Harmonic Analysis},
  volume  = {42},
  number  = {3},
  pages   = {378--401},
  year    = {2017},
  doi     = {10.1016/j.acha.2015.08.014}
}

@article{ADK2013,
  author  = {Aldroubi, Akram and Davis, Jacqueline and Krishtal, Ilya},
  title   = {Dynamical sampling: time--space trade-off},
  journal = {Applied and Computational Harmonic Analysis},
  volume  = {34},
  pages   = {495--503},
  year    = {2013},
  number        = {3},
  doi           = {10.1016/j.acha.2012.09.002},
}

@article{ADK2015,
  author  = {Aldroubi, Akram and Davis, Jacqueline and Krishtal, Ilya},
  title   = {Exact reconstruction of signals in evolutionary systems via
             spatiotemporal trade-off},
  journal = {Journal of Fourier Analysis and Applications},
  volume  = {21},
  pages   = {11--31},
  year    = {2015},
  number        = {1},
  doi           = {10.1007/s00041-014-9359-9},
}

@article{CMPP2020,
  author  = {Cabrelli, Carlos and Molter, Ursula and Paternostro, Victoria and
             Philipp, Friedrich},
  title   = {Dynamical sampling on finite index sets},
  journal = {Journal d'Analyse Math\'ematique},
  volume  = {140},
  pages   = {637--667},
  year    = {2020},
  doi           = {10.1007/s11854-020-0099-2},
}

@misc{ACNP2025,
  author        = {Aguilera, Alejandra and Cabrelli, Carlos and Negreira, Felipe and
                   Paternostro, Victoria},
  title         = {Optimal dynamical frames},
  year          = {2025},
  eprint        = {2506.00567},
  archivePrefix = {arXiv},
  primaryClass  = {math.FA},
  note          = {},
  doi           = {},
}

@ARTICLE{AP2017,
       author = {{Aldroubi}, Akram and {Petrosyan}, Armenak},
        title = "{Dynamical sampling and systems from iterative actions of operators}",
      journal = {arXiv e-prints},
         year = 2016,
          eid = {arXiv:1606.03136},
        archivePrefix = {arXiv},
       eprint = {1606.03136},
        primaryClass = {math.FA},
       adsurl = {https://ui.adsabs.harvard.edu/abs/2016arXiv160603136A}
}

@article{ACCMP2017normal,
  author  = {Aldroubi, Akram and Cabrelli, Carlos and \c{C}akmak, Ahmet Faruk and
             Molter, Ursula and Petrosyan, Armenak},
  title   = {Iterative actions of normal operators},
  journal = {Journal of Functional Analysis},
  volume  = {272},
  pages   = {1121--1146},
  year    = {2017},
  number        = {3},
  doi           = {10.1016/j.jfa.2016.10.027},
}

@article{CHP2020,
  author  = {Christensen, Ole and Hasannasab, Marzieh and Philipp, Friedrich},
  title   = {Frame properties of operator orbits},
  journal = {Mathematische Nachrichten},
  volume  = {293},
  number  = {1},
  pages   = {52--66},
  year    = {2019},
  doi     = {10.1002/mana.201800344},
}

@inbook{CMS2023,
author = {Cabrelli, Carlos and Molter, Ursula and Suárez, Daniel},
year = {2024},
month = {01},
pages = {329-348},
title = {Frames of Iterations and Vector-Valued Model Spaces},
isbn = {978-3-031-41129-8},
doi = {10.1007/978-3-031-41130-4_12}
}

@misc{ADGMM2023,
  author        = {Aldroubi, Akram and D\'iaz Mart\'in, Roc\'io and Gong, Le and
                   Mashreghi, Javad and Medri, Ivan},
  title         = {Characterization of frames for source recovery from dynamical
                   samples},
  year          = {2024},
  eprint        = {2401.15450},
  archivePrefix = {arXiv},
  primaryClass  = {math.DS},
  doi           = {10.48550/arXiv.2401.15450},
}

@article{ACM2025periodic,
  author  = {Aldroubi, Akram and Cabrelli, Carlos and Molter, Ursula},
  title   = {Periodic source detection in discrete dynamical systems via
             space--time sampling},
  journal = {Expositiones Mathematicae},
  volume  = {43},
  pages   = {Paper No. 125730},
  year    = {2025},
  number        = {6},
  doi           = {10.1016/j.exmath.2025.125730},
}

@article{ADMM2024,
  author  = {Aldroubi, Akram and D\'iaz Mart\'in, Roc\'io and Medri, Ivan},
  title   = {Dynamical sampling for the recovery of spatially constant source
             terms in dynamical systems},
  journal = {Linear Algebra and its Applications.},
  volume  = {694},
  pages   = {148--185},
  year    = {2024},
  doi           = {10.1016/j.laa.2024.04.011},
}

@article{AHKK2023burst,
  author  = {Aldroubi, Akram and Huang, Longxiu and Kornelson, Keri and Krishtal, Ilya},
  title   = {Predictive algorithms in dynamical sampling for burst-like forcing
             terms},
  journal = {Applied and Computational Harmonic Analysis},
  volume  = {65},
  pages   = {322--347},
  year    = {2023},
  doi           = {10.1016/j.acha.2023.03.003},
}

@incollection{DMM2021control,
  author="Mart{\'i}n, Roc{\'i}o D{\'i}az
and Medri, Ivan
and Molter, Ursula",
editor="Hirn, Matthew
and Li, Shidong
and Okoudjou, Kasso A.
and Saliani, Sandra
and Yilmaz, {\"O}zg{\"u}r",
title="Dynamical Sampling: A View from Control Theory",
bookTitle="Excursions in Harmonic Analysis, Volume 6: In Honor of John Benedetto's 80th Birthday",
year="2021",
publisher="Springer International Publishing",
address="Cham",
pages="269--295",
isbn="978-3-030-69637-5",
doi="10.1007/978-3-030-69637-5_15",
url="https://doi.org/10.1007/978-3-030-69637-5-15"
}

@article{AHKLLV2018noise,
  author  = {Aldroubi, Akram and Huang, Longxiu and Krishtal, Ilya and Ledeczi, Akos and Lederman, Roy R. and Volgyesi, Peter},
  title   = {Dynamical sampling with additive random noise},
  journal = {Sampling Theory in Signal and Image Processing},
  volume  = {17},
  pages   = {153--182},
  year    = {2018},
  doi     = {10.1007/BF03549662}
}

@article{Deutsch1995,
author = {Deutsch, Frank},
year = {1995},
month = {01},
pages = {},
title = {The Angle Between Subspaces of a Hilbert Space},
isbn = {978-90-481-4516-4},
journal = {Journal of Approximation Theory - JAT},
doi = {10.1007/978-94-015-8577-4_7}
}

@article{Eldar2003,
  author  = {Eldar, Yonina C.},
  title   = {Sampling with arbitrary sampling and reconstruction spaces and
             oblique dual frame vectors},
  journal = {Journal of Fourier Analysis and Applications},
  volume  = {9},
  number  = {1},
  pages   = {77--96},
  year    = {2003},
  doi           = {10.1007/s00041-003-0004-2},
}

@article{Stuart2010,
  author  = {Stuart, Andrew M.},
  title   = {Inverse problems: a {B}ayesian perspective},
  journal = {Acta Numer.},
  volume  = {19},
  pages   = {451--559},
  year    = {2010},
  doi           = {10.1017/S0962492910000061},
}

@book{RW2006,
  author    = {Rasmussen, Carl Edward and Williams, Christopher K. I.},
  title     = {Gaussian Processes for Machine Learning},
  publisher = {MIT Press},
  year      = {2006},
  doi           = {10.7551/mitpress/3206.001.0001},
}

@article{ALL2013,
  author  = {\'Alvarez, Mauricio A. and Luengo, David and Lawrence, Neil D.},
  title   = {Linear latent force models using {G}aussian processes},
  journal = {IEEE Transaction on Pattern Analysis and Machine Intelligence},
  volume  = {35},
  number  = {11},
  pages   = {2693--2705},
  year    = {2013},
  doi           = {10.1109/TPAMI.2013.86},
}

@book{AdlerTaylor2007,
  author    = {Adler, Robert J. and Taylor, Jonathan E.},
  title     = {Random Fields and Geometry},
  publisher = {Springer},
  year      = {2007},
  doi           = {10.1007/978-0-387-48116-6},
}

@article{AK2016krylov,
  author  = {Aldroubi, Akram and Krishtal, Ilya},
  title   = {Krylov subspace methods in dynamical sampling},
  journal = {Sampling Theory in Signal and Image Processing},
  volume  = {15},
  pages   = {9--20},
  year    = {2016},
  doi           = {10.1007/BF03549595},
}

@article{Tang2017sysid,
  author  = {Tang, Sui},
  title   = {System identification in dynamical sampling},
  journal = {Advances in Computational Mathematics},
  volume  = {43},
  pages   = {555--580},
  year    = {2017},
  number        = {3},
  doi           = {10.1007/s10444-016-9497-5},
}

@article{KMT2024learning,
  author  = {K\"ummerle, Christian and Maggioni, Mauro and Tang, Sui},
  title   = {Learning transition operators from sparse space-time samples},
  journal = {IEEE Transaction on Information Theory},
  volume  = {70},
  pages   = {6412--6446},
  year    = {2024},
  number        = {9},
  doi           = {10.1109/TIT.2024.3413534},
}

@article{DMO2023error,
  author  = {D\'iaz Mart\'in, Roc\'io and Medri, Ivan and Osorio, Juliana},
  title   = {Error analysis on the initial state reconstruction problem},
  journal = {Sampling Theory, Signal Processing and Data Analysis},
  volume  = {21},
  pages   = {Paper No. 18},
  year    = {2023},
  number        = {1},
  doi           = {10.1007/s43670-023-00059-5},
}

\end{document}